\documentclass{amsart}
\usepackage{amsxtra,amssymb}
\theoremstyle{plain}

\newcommand{\cleqn}{\setcounter{equation}{0}}
\newcommand{\clth}{\setcounter{theorem}{0}}
\newcommand {\sectionnew}[1]{\section{#1}\cleqn\clth}
\newcommand{\nn}{\hfill\nonumber}
\newtheorem{theorem}{Theorem}[section]
\newtheorem{lemma}[theorem]{Lemma}
\newtheorem{definition-theorem}[theorem]{Definition-Theorem}
\newtheorem{proposition}[theorem]{Proposition}
\newtheorem{corollary}[theorem]{Corollary}
\newtheorem{definition}[theorem]{Definition}
\newtheorem{example}[theorem]{Example}
\newtheorem{remark}[theorem]{Remark}
\newcommand \bth[1] { \begin{theorem}\label{t#1} }
\newcommand \ble[1] { \begin{lemma}\label{l#1} }

\newcommand \bpr[1] { \begin{proposition}\label{p#1} }
\newcommand \bco[1] { \begin{corollary}\label{c#1} }
\newcommand \bde[1] { \begin{definition}\label{d#1}\rm }
\newcommand \bex[1] { \begin{example}\label{e#1}\rm }
\newcommand \bre[1] { \begin{remark}\label{r#1}\rm }

\renewcommand {\eth} { \end{theorem} }
\newcommand {\ele} { \end{lemma} }

\newcommand {\epr} { \end{proposition} }
\newcommand {\eco} { \end{corollary} }
\newcommand {\ede} { \end{definition} }
\newcommand {\eex} { \end{example} }
\newcommand {\ere} { \end{remark} }

\newcommand \thref[1]{Theorem \ref{t#1}}
\newcommand \leref[1]{Lemma \ref{l#1}}
\newcommand \prref[1]{Proposition \ref{p#1}}
\newcommand \coref[1]{Corollary \ref{c#1}}
\newcommand \deref[1]{Definition \ref{d#1}}

\newcommand \lb[1]{\label{#1}}
\def \Rset {{\mathbb R}}         
\def \Cset {{\mathbb C}}
\def \Zset {{\mathbb Z}}

\def \Qset {{\mathbb Q}}
\def \Pset {{\mathbb P}}
\def \cp   {{\Cset \Pset^1}}
\def \g  {\mathfrak{g}}   
\def \gt  {\tilde{\mathfrak{g}}}   
\def \k  {\mathfrak{k}}
\def \p  {\mathfrak{p}}  
\def \pa  {\tilde{\mathfrak{p}}}  
\def \ja  {\hat{j}}         
\def \b  {\mathfrak{b}}  
\def \q  {\mathfrak{q}}  
\def \f  {\mathfrak{f}}  
\def \n  {\mathfrak{n}}  
\def \na  {\tilde{\mathfrak{n}}}  
\def \m  {\mathfrak{m}}  
\def \a  {\mathfrak{a}}
\def \h  {\mathfrak{h}}
\def \t  {\mathfrak{t}}

\def \C  {{\mathcal{C}}}           
\newcommand \Aff  {\mathcal{AFF}}
\def \Fin  {\mathcal{F} in}

\def \H  {{\mathcal{H}}}
\def \GK {{\mathcal{C}_{(\g, \k)}}}

\def \GtK {{\mathcal{C}_{(\gt, \k, \kappa)}}}
\def \PK  {{\mathcal{C}_{(\pa_+, \k, \kappa)}}}
\def \KL  {{\mathcal{KL}}}

\def \O  {{\mathcal{O}}}

\def \c  {{K}}

\def \De {\Delta}   
\def \hchi {\widehat{\chi}}
\def \de {\delta}
\def \al {\alpha}
\def \be {\beta}
\def \Ga {\Gamma}
\def \ga {\gamma}
\def \ka {\kappa}
\def \Om {\Omega}
\def \om {\omega}
\def \e  {t}   
\def \la {\lambda}
\def \La {\Lambda}

\def \vph {\varphi}
\def \Sig {\Sigma}

\def \ig {inf}
\def \int {int}

\def \mt  {\mapsto}

\def \ra  {\rightarrow}           
           
\def \thra {\twoheadrightarrow}

\def \hra {\hookrightarrow}
\def \sub {\subset}

\def \all {\langle}
\def \arr {\rangle}
\def \st  {\ast}                 

\def \wh {\widehat}

\def \id { {\mathrm{id}} }
\def \Ann { {\mathrm{Ann}} }

\DeclareMathOperator \Span { {\mathrm{span}} }

\DeclareMathOperator \Hom { {\mathrm{Hom}} }
\DeclareMathOperator \Ext { {\mathrm{Ext}} }
\DeclareMathOperator \I { {\mathrm{I}} }
\DeclareMathOperator \Ind { {\mathrm{Ind}} }
\DeclareMathOperator \Irr { {\mathrm{Irr}} }
\DeclareMathOperator \Res { {\mathrm{Res}} }

\renewcommand \Re { {\mathrm{Re}} }

\newcommand \ev { {\mathrm{Ex}} }
\begin{document}
\title[Affine Jacquet functors]
{
Affine Jacquet functors and \\ Harish-Chandra categories
}
\author[Milen Yakimov]{Milen Yakimov}
\address{
Department of Mathematics \\
University of California \\
Santa Barbara, CA 93106, U.S.A.}
\email{yakimov@math.ucsb.edu}
\date{}
\begin{abstract}
We define an affine Jacquet functor and use it to describe the
structure of induced affine Harish-Chandra modules
at noncritical levels, 
extending the theorem of Kac and Kazhdan \cite{KK} on 
the structure of Verma modules in the Bernstein--Gelfand--Gelfand 
categories $\O$ for Kac--Moody algebras. This is combined with
a vanishing result for certain extension groups to construct 
a block decomposition of the categories of affine Harish-Chandra 
modules of Lian and Zuckerman \cite{LZ}. The latter provides 
an extension of the works of Rocha-Caridi, Wallach \cite{RW} and 
Deodhar, Gabber, Kac \cite{DGK} on block decompositions of
BGG categories for Kac-Moody algebras. We also prove a compatibility 
relation between the affine Jacquet functor and the Kazhdan--Lusztig 
tensor product. A modification of this is used to prove 
that the affine Harish-Chandra category is stable under 
fusion tensoring with the Kazhdan--Lusztig category 
(a case of our finiteness result \cite{Y}) and will be 
further applied in studying translation functors
for Kac--Moody algebras, based on the fusion tensor 
product.
\end{abstract}
\maketitle
\sectionnew{Introduction}\lb{Intro}
Let $\g$ be a complex semisimple Lie algebra and $\gt$ be the 
corresponding affine Kac--Moody algebra.

In this paper we define an affine version of the 
Jacquet functor \cite{Cas, W}, introduced by Casselman and Wallach, 
and use it to deduce properties of the affine Harish-Chandra categories
and the Kazhdan--Lusztig fusion tensor product.

The standard Harish-Chandra category, associated to 
a real form $\g_0$ of $\g$ (consisting of finite length,
admissible $(\g, \k)$-modules) will be denoted by $\H.$
The Jacquet module of $V \in H$ is given by
\[
j(V) = \lim_{\ra} \Ann_{\n_0^k} V^\st \sub V^\st
\]
where $\n_0$ is the nil radical of a minimal parabolic 
subalgebra $\q_0$ of $\g_0.$ It is a faithful and exact
(contravariant) functor from $\H$ to a generalized
Bernstein--Gelfand--Gelfand category, to be denoted
by $\O'.$ (On the latter $(\q_0/\n_0)_\Cset$ only
acts locally finitely, but in general not semisimply.)

The affine Harish-Chandra category $\wh{\H}_\ka$ consists of 
finitely generated
smooth $\gt$-modules of central charge $\ka - h\spcheck$ on which 
the Sugawara operator $L_0$ acts locally finitely and the
corresponding generalized eigenspace decomposition of $L_0$
has the form
\begin{equation}
\label{HHH}
V = \bigoplus_{\xi \colon \xi - \xi_i \in \Zset_{\geq 0} }
V^\xi \; \mbox{for some} \; \xi_1, \ldots, \xi_n \in \Cset
\end{equation}
where $V^\xi \in \H,$ considered as $\g$-modules.
These categories were introduced by Lian and Zuckerman
\cite{LZ, LZ2} in a slightly different way. See section 2.3
for details. (As usual $h\spcheck$ denotes the dual Coxeter
number of $\g.)$

The generalized affine Bernstein--Gelfand--Gelfand category
$\wh{\O}'_\ka$ is the category of smooth $\gt$-modules
of central charge $\ka - h\spcheck$ with the 
above property \eqref{HHH} and $V^\xi \in \O'.$

The affine Jacquet functor is a faithful and exact 
(contravariant) functor 
$\ja \colon \wh{\H}_\ka \ra \wh{\O}'_\ka,$ given by
\begin{equation}
\label{jaja}
\ja(V) = \left[ j(V)^\sharp(\infty) \right]^{\Cset[L_0]-fin}.
\end{equation}
Here $(.)^\sharp$ denotes the twist of a $\gt$-module
by the automorphism of $\g:$ 
$(x t^n)^\sharp = x (-t)^n,$ $x \in \g,$ $K^\sharp = - K.$
By $(.)(\infty)$ we denote the strictly smooth part 
of a $\gt$-module, see section 2.3. The notation 
$W^{\Cset[L_0]-fin}$ is for the $\Cset[L_0]$ locally 
finite part of $W.$ {\em{This can be omitted in the 
case $\ka \notin \Qset_{\geq 0}.$}}

If the $\gt$-module $V \in \wh{\H}_\ka$ has the $L_0$ 
generalized eigenspace decomposition \eqref{HHH}, then
\[
\ja(V) = \bigoplus_{\xi \colon \xi - \xi_i \in \Zset_{\geq 0} } 
\lim_{\ra} \Ann_{\n_0^k} (V^\xi)^\st
\]
as $\gt$-submodules of $(V^\st)^\sharp$ where $(V^\xi)^\st$ is identified 
with the subspace of $V^\st,$ consisting of functionals vanishing on 
$\bigoplus_{\xi' \neq \xi} V^{\xi'}.$

For every irreducible $\g$-module $M$ the induced
$\gt$-module
\[
\Ind(M)= U(\gt) \otimes_{U(\g[t] \oplus \Cset K)} M
\]
and its unique irreducible quotient $\Irr(M)_\ka$ belong
to $\wh{\H}_\ka.$ Here $K$ acts on $M$ by 
$\ka -h\spcheck$ and $\g[t]$ acts on $M$ through the quotient
$\g[t] \thra \g[t]/t \g[t] \cong \g.$

We prove that
\begin{equation}
\label{iii}
\ja( \Ind(M)_\ka ) = D(\Ind(j(M)^d)_\ka)
\end{equation}
where $(.)^d$ and $D(.)$ refer to certain natural duality
functors in $\H$ and $\wh{\H}_\ka,$ respectively, see
section 2.1 and 2.3.

The idea to use the Jacquet functor to study 
the structure of induced modules 
in the affine Harish-Chandra categories
belongs to Lian and Zuckerman \cite{LZ2}. They formulated
a version of \eqref{iii} but with the use of 
the standard Jacquet functor
which makes it slightly incorrect. 

{} From \eqref{iii} we prove a generalization to the affine 
Harish-Chandra categories of
the result of Kac and Kazhdan on the structure of Verma modules
for Kac--Moody algebras:

{\em{If $M$ is an irreducible Harish-Chandra module with 
infinitesimal character $\chi_\la \in \h^\st/W$ 
for some $\la \in \h^\st,$ then all irreducible subquotients  
of $\Ind(M)_\ka$ are isomorphic to $\Ind(M_2)_\ka$
for some irreducible Harish-Chandra modules $M_2$
with infinitesimal character $\chi_{\la_2}$
such that there exists $w \in W$ for which
the pair $[w \la, \la_2] \in \h^\st \times \h^\st$ 
satisfies the condition
$(\st)$ of Kac-Kazhdan, reviewed in \deref{star}.
}}

We further obtain in section 4.3 a block decomposition of the categories
$\wh{\H}_\ka,$ extending the works of Rocha-Caridi, Wallach \cite{RW} and 
Deodhad, Gabber, Kac \cite{DGK} on block decomposition of the categories 
$\O$ for Kac--Moody algebras.
{\em{Define the equivalence relation $\sim$ in $\h,$
induced by $\la \sim \mu$ if $\mu \in W \la$ or the pair 
$[\la, \mu] \in \h^\st \times \h^\st$ satisfies the condition 
$(\st)$ in \deref{star}. Then the categories of affine 
Harish-Chandra modules $\wh{\H}_\ka$ posses the block decompositions
\[
\wh{\H}_\ka = \bigoplus_{\hchi \in (\h^\st/\sim)} 
\wh{\H}_\ka^{\hchi}.
\]
The subcategories $\wh{\H}_\ka^{\hchi}$ consist
of modules $V \in \wh{\H}_\ka$ with a filtration by
$\gt$-submodules
\[
0=W_0 \sub W_1 \sub \ldots \sub W_N =V
\]
for which the subquotients $W_i/W_{i-1}$ are isomorphic
to quotients of $\Ind(M)_\ka$ for irreducible Harish-Chandra
modules with infinitesimal characters in the class
$\hchi \in (\h^\st/\sim).$}} Similar result is proved for
infinitely generated affine Harish-Chandra modules in analogy 
with \cite{DGK}.

{\em{We also show vanishing of $\Ext$ groups between different 
blocks in a much larger category}} 
in the spirit of the results of 
Rocha-Caridi and Wallach \cite{RW}, and the vanishing
of $\Ext$'s between blocks of $\H$ in the larger category
of all $(\g, \k)$-modules (neither finitely generated, nor
admissible, in general), see \cite[Theorem 4.1, Ch 1]{BW}
and \prref{Extvan} below. This is done in section 4.2.

Our approach is very similar to the one of Rocha-Caridi,
Wallach \cite{RW} and Deodhar, Gabber, Kac \cite{DGK}. 
(The paper \cite{DGK} discusses only the case
of vanishing of $\Ext^1$ groups in the the
$\O$ category for an arbitrary Kac--Moody algebra 
and \cite{RW} the general case.) The proof of the vanishing of 
$\Ext$ groups in the category $\O$ in \cite{RW, DGK} 
uses the fact that modules from this category, 
restricted to the Cartan subalgebra of the extended 
affine Kac--Moody algebra belong to a semisimple category.
On the contrary modules from the category $\wh{\H}_\ka,$ 
which we treat, restricted to $\g \hra \hat{\g}$ form essentially the 
Harish-Chandra category (which is non-semisimple). The use of 
semisimplicity in \cite{DGK, RW} can be avoided. When written in terms  
of spectral sequences the arguments of \cite{RW, DGK} simplify a lot,
as it is always the case in similar situations.

In section 5 we derive the following compatibility property
between the affine Jacquet functor and the Kazhdan--Lusztig
tensor product:

{\em{Let $\ka \notin \Qset_{\geq 0}.$ For any module $U$ 
in the Kazhdan--Lusztig category and $V \in \wh{\H}_\ka$
the following isomorphism holds:}}
\begin{equation}
\label{jjj}
\ja(U \dot{\otimes} V ) \cong
D \left[ U \dot{\otimes} D \ja(V) \right].
\end{equation}

In \cite{Y} we showed that for any subalgebra $\f$ of $\g$ which is
reductive in $\g$ the affine analogs of the categories of finite length,
admissible $(\g, \f)$-modules are stable under the fusion tensoring
with the Kazhdan--Lusztig category for $\ka \notin \Rset_{\geq 0}.$
This is an affine version Kostant's theorem \cite{Kos}.

A modification of \eqref{jjj} provides another proof of a 
case of our result \cite{Y}, namely that the affine Harish-Chandra
category is stable under fusion tensoring with modules from 
the Kazhdan--Lusztig category. 

The above compatibility of the fusion tensor product and the affine 
Jacquet functor will be further used in studying fusion translation
functors in the spirit of Zuckerman \cite{Z} and Jantzen \cite{J}.

After the work of Beilinson and Bernstein \cite{BB} a very general 
nonvanishing of the Jacquet modules of a $\g$-module with an infinitesimal 
character (associated to the nilradicals of Borel subalgebras 
in a dense subset of 
the flag variety) is known. It is interesting to understand if this
can be used to define block decomposition of the category of all
smooth $(L_0$-locally finite) $\gt$-modules for which $V^\xi$ 
are finitely generated $U(\g)$-modules on which the center 
of $U(\g)$ acts locally finitely. In this way the affine Jacquet
functor would substitute very efficiently the triviality of the 
center of an affine Kac--Moody algebra.

Acknowledgments: I am indebted to Dan Barbasch, Yuri Berest, and
Edward Frenkel for very helpful discussions and comments.
\sectionnew{Preliminaries on real semisimple Lie algebras and affine 
Kac-Moody algebras}
\subsection{Real semisimple Lie algebras}
Let $\g_0$ be a real semisimple Lie algebra with complexification
$\g = (\g_0)_\Cset.$ Fix a Cartan decomposition of $\g_0$
\[
\g_0 = \k_0 \oplus \p_0.
\]
Let $\a_0 \sub \p_0$ be a maximal commutative subalgebra. Denote by
$\Sig_+$ a fixed set of positive restricted roots of $\g_0$ with respect 
to $\a_0$ and set
\[
\n_0 = \bigoplus_{\la \in \Sig_+} \g_0^\la.
\]
Define the minimal parabolic subalgebra of $\g_0$ 
\[
\q_0 = \m_0 \oplus \a_0 \oplus \n_0
\]
where $\m_0 = Z_{\k_0}(\a_0).$ 

Fix a maximal commutative subalgebra $\t_0 \sub \m_0$ and consider the 
related maximally noncompact Cartan subalgebra 
\[
\h_0 = \t_0 \oplus \a_0
\]
of $\g_0.$ Then the complexification $\h= (\h_0)_\Cset$ is 
a Cartan subalgebra of $\g$ and one can find a set of positive roots
$\Delta_+$ for $\g$ relative to $\h$ which extends $\Sig_+,$ i.e.
\[
\Delta|_{\a_0} = \Sig_+ \cup \{ 0 \}.
\]
Denote the nil radical
\[
\n_+ = \bigoplus_{\la \in \Delta_+} \g^\la
\]
of the related positive Borel subalgebra. Then
\[
\n_+ \supset (\n_0)_\Cset.
\]

Recall that the category $\H$ of Harish-Chandra modules for the 
pair $(\g, \k)$ is the category of finitely generated $\g$-modules
$V$ such that
\[
V|_{\k} = \bigoplus_{\mu} (V^\mu)^{\oplus m_\mu}
\]
where $V^\mu$ are representatives of the equivalence classes
of all irreducible finite dimensional $\k$-modules  
and the multiplicities $m_\mu$ are finite (admissibility 
condition). All Harish-Chandra modules 
have finite length.

The generalized Bernstein--Gelfand--Gelfand category
$\O'$ for $\g$ related to the choice of Borel subalgebra $\b$ and 
Cartan subalgebra $\h$ as above, consists of finitely generated 
$\g$-modules $V$ on which $\n_+$ acts locally nilpotently and 
\[
V|_{\h} = \bigoplus_{\la \in \h^*} V_\la, \quad \dim V_\la < \infty
\]
where $V_\la$ denote the generalized eigenspaces of $\h:$
\[
V_\la = \{ v \in V \mid (h - \la(h))^k v = 0, \, \forall h \in \h \; 
\mbox{for some} \; k \in \Zset_{>0} \}.
\]
The usual BGG category for $\g$ for the choices of Cartan and Borel 
subalgebras made will be denoted by $\O.$

The center of $\g$ will be denoted by $Z(\g).$ By the Harish-Chandra
isomorphism $Z(\g) \cong S(\h)^W$ and the characters of  
$Z(\g)$ are parametrized by $\h^\st/W.$ The character corresponding
to $\la \in \h^\st$ will be denoted by $\chi_\la \colon Z(\g) \ra \Cset.$
Recall that for the Casimir element of $U(\g)$
\begin{equation}
\label{chiOm}
\chi_\la(\Om) = |\la|^2 - |\rho|^2.
\end{equation}
where $\rho$ is the half-sum of the positive roots of the Borel subalgebra 
of $\g$ used to define the Harish-Chandra isomorphism.

For a $\g$-module $V$ possessing an infinitesimal character, 
by $\chi(V)$ we will denote the latter.

The categories $\H$ and $\O'$ posses block decompositions
\begin{equation}
\label{block}
\H = \bigoplus_{\chi \in \h^\st/W} \H^\chi, \quad
\O' = \bigoplus_{\chi \in \h^\st/W} \O^{' \chi}.
\end{equation}
where $\H^\chi$ and $\O^{' \chi}$ are the full subcategories
of $\H$ and $\O',$ respectively, consisting of modules with
filtrations by $\g$-modules with infinitesimal character
$\chi.$

Later we will need an important result on vanishing of $\Ext$ groups
between Harish-Chandra modules modules from different blocks of $\H,$
see \cite[Chapter I, Theorem 4.1]{BW}. Denote by $\GK$ the category of 
$(\g, \k)$ modules, i.e. $\g$-modules which are locally $\k$-finite and 
$\k$-semisimple. It is well known that these categories have enough 
projectives and injectives, see e.g. \cite[Chapter I, 2.4]{BW}.

\bpr{Extvan} Assume that $V_i \in \H^{\chi_i},$ $i=1,2$ and 
$\chi_1 \neq \chi_2.$ Then
\[
\Ext^n_{\g, \k}(V_1, V_2) =0, \; n \in \Zset_{\geq 0}
\]
where $\Ext^n_{\g, \k}$ refer to the $\Ext$ groups in the category
$\GK.$
\epr

The categories $\H$ and $\O'$ have natural duality functors 
(involutive antiequivalences). Both of them will be denoted by
$V \mt V^d.$ They are given by
\begin{align}
V^d &= (V^\st)^{U(\k)-fin}, \quad V \in \H,
\label{dualH}
\\
V^d &= (V^\st)^{U(\h)-fin}, \quad V \in \H.
\label{dualO}
\end{align}
Here and later for any module $W$ over a Lie algebra by $W^\st$ we denote 
the (full) dual module. For a $\Cset$-algebra $A$ and an $A$-module
$W$ we denote by $W^{A-fin}$ the $A$-submodule of $W$ consisting
of $A$-finite vectors, i.e.
\[
W^{A-fin} := \{ w \in W \mid \dim A w < \infty \}.
\]

Clearly $(.)^d$ restricts to a functor from 
$\H^{\chi_\la}$ to $\H^{\chi_{-\la}}$ and from 
$\O^{' \chi_\la}$ to $\O^{' \chi_{-\la}}.$
\subsection{The Jacquet functor}
Fix $V \in \H.$ The natural increasing subsequence of 
$\m_0 \oplus \a_0$-submodules of $V$
\[
\n_0 V \supset \n_0^2 V \supset \ldots
\]
gives rise to the increasing subsequence of 
$\m_0 \oplus \a_0$-submodules of $V^\st$ 
\[
(V/\n_0 V)^\st \hra  (V/\n_0^2 V)^\st \hra \ldots .
\]
Here for a subspace $V_1$ of $V,$ $(V/V_1)^\st$ is naturally identified 
with the subspace of $V^\st$ that consists of $\eta \in V^\st$ such that 
$\eta|_{V_1}=0.$ Note that
\[
(V/\n_0^kV)^\st = \Ann_{\n_0^k} V^\st.
\]
Finally the Jacquet module of $V$ is defined by
\[
j(V) = \lim_{\ra} \Ann_{\n_0^k} V^\st = 
\lim_{\ra} (V/\n_0^k V)^\st.
\]
Since the adjoint action of $\n_0$ on $\g_0$ is nilpotent 
the space $j(V)$ is a $\g$-submodule of $V^\st.$

\bth{Jacquet} (i) For any $V \in \H,$ 
$\dim (V/\n_0^k V) < \infty$ and $j(V) \in \O'.$ 

(ii) The contravariant functor $j \colon \H \ra \O'$ is 
faithful and exact.
\eth

Because of part (i) of the above theorem, the Jacquet functor
is also given by 
\[
j(V) = \lim_{\ra} \Ann_{\n_+^k} V^\st = 
\lim_{\ra} (V/\n_+^k V)^\st,
\]
see section 2.1 for the definition of $\n_+.$

In fact the Jacquet functor takes values in the subcategory 
of $\O'$ consisting of $\g$-modules on which the Levi factor
$(\m_0 \oplus \a_0)_\Cset$ of the parabolic subalgebra 
$\q= (\q_0)_\Cset$ of $\g$ acts locally finitely. 

It is also clear that 
$j \colon \H^{\chi_\la} \ra {\O}^{' \chi_{-\la}}.$ 

Recall that the Harish-Chandra category $\H$ 
is stable under tensoring with finite dimensional modules.
We have the following property, relating this tensoring
with the Jacquet functor.

\bpr{finJa} For any finite dimensional $\g$-module $U$ and a 
Harish-Chandra module $V$ the Jacquet module of $U \otimes V \in \H$ 
is given by
\[
j(U \otimes V) \cong U^\st \otimes j(V).
\]
\epr

\subsection{Untwisted Affine Kac--Moody algebras}
Consider the loop algebra 
$\g [t, t^{-1}] = \g \otimes_\Cset \Cset[t, t^{-1}].$
The untwisted affine Kac-Moody algebra $\gt$ associated to $\g$ is the
central extension of this loop algebra by
\begin{equation}
\label{affine}
[x t^n, y t^m] = [x, y] t^{n+m} + n \de_{n, -m} (x, y) \c, \;
x, y \in \g
\end{equation}
where $(., .)$ is an invariant nondegenerate bilinear form on $\g$
normalized by $(\al, \al)=2$ for long roots $\al$ of $(\g, \h).$
For a full exposition of Kac--Moody algebras we refer to Kac's
book \cite{Kac}.

A $\gt$-module $V$ is called smooth if any $v \in V$ is annihilated
by $x t^n$ for all $x \in \g$ and $n \gg 0.$ 

On any smooth $\gt$-module of central charge
$\ka - h\spcheck$ $(\ka \neq 0)$ one has a natural representation
of the Virasoro algebra \cite{Kac}, given by the Sugawara 
operators
\begin{equation}
\label{Sug}
L_k = \frac{1}{2 \ka} \sum_p \sum_{j \in \Zset}
\colon (x_p t^{-j}) (x_p t^{j+k}) \colon.
\end{equation}
Here $\{x_p\}$ is an orthonormal basis of $\g$ 
with respect to the bilinear form $(.,.).$ 
The normal ordering in \eqref{Sug} prescribes pulling to the right 
the term $x t^n$ with larger n. Here and later $h\spcheck$ denotes the 
dual Coxeter number of $\g.$ The following commutation relations hold
\begin{equation}
\label{commut}
[L_k, x t^n] = -n (x t^{n+k}).
\end{equation}
For any smooth $\gt$-module $V$ define the generalized eigenspaces 
of $L_0$
\[
V^\xi = \{ v \in V \mid (L_0 - \xi)^n v = 0, \; 
\mbox{for some integer} \, n \}, \, \xi \in \Cset.
\]
Because $\g \hra \gt$ commutes with $L_0$ every $V^\xi \sub V$ is 
naturally a $\g$-module.

\bde{Aff} Let $\C$ be a category of finitely generated $\g$-modules that 
is closed under tensoring with finite dimensional $\g$-modules. 
Define the category $\Aff(\C)_\ka$ to be the representation category of 
finitely generated smooth $\gt$-modules of central charge $\ka - h\spcheck$
which are $L_0$ locally finite and satisfy
\begin{equation}
\label{Sugdecom}
V = \bigoplus_{\xi \colon \xi-\xi_i \in \Zset_+} V^\xi \quad
\mbox{for some} \quad \xi_1, \xi_2, \ldots, \xi_n \in \Cset
\end{equation}
with $V^\xi \in \C,$ considered as $\g$-modules.
\ede 

If $\Fin$ is the category of finite dimensional $\g$-modules then the 
Kazhdan--Lusztig categories \cite{KL1} for $\gt$ are 
$\KL_\ka = \Aff(\Fin)_\ka$ in the case $\ka \notin \Qset_{\geq 0}.$

Denote 
\[
\wh{\O}_\ka = \Aff(\O)_\ka, \quad  
\wh{\O}'_\ka = \Aff(\O')_\ka, \quad 
\ka \neq 0.
\]  

Recall that the extended Kac-Moody algebra $\hat{\g}$ associated to $\g$
is the Lie algebra $\gt \oplus \Cset d$ where 
\[
[d, x t^n] = n x t^n, \quad [d, K] =0,
\] 
see \cite{Kac} for details.
Consider its Cartan subalgebra 
$\hat{\h}=\h \oplus \Cset K \oplus \Cset d$ and 
Borel subalgebra $\b_+ \oplus \g[t] \oplus \Cset K \oplus \Cset d.$  
The associated affine BGG category of noncritical central charge 
$\ka - h\spcheck$ is essentially the category $\wh{O}_\ka$ with the 
difference that the generator $d$ 
can act on any module $V \in \wh{O}_\ka$ 
by $const - L_0$ for any choice of the constant involved. 
  
The categories of affine Harish-Chandra modules defined by Lian and 
Zuckerman \cite{LZ, LZ2} are $\Aff(\H)_\ka.$ They will be denoted
by $\wh{\H}_\ka.$

The Lie algebra $\gt$ is $\Zset$-graded by
\begin{equation}
\label{grading}
\deg x t^n = -n, n \in \Zset, x \in \g; \;
\deg K =0.
\end{equation}
Each module in the categories $\Aff(\C)_\ka$ is naturally 
$\Cset$-graded with respect to the grading \eqref{grading} 
of $\gt$ using the generalized eigenvalue decomposition 
\eqref{Sugdecom} of $L_0.$ This is the reason for the choice of the 
negative sign in \eqref{grading}. Moreover each morphism 
in $\Aff(\C)_\ka$ preserves the grading \eqref{Sugdecom}. 

For simplicity we will denote the maximal parabolic subalgebra of $\gt$
\begin{equation}
\pa_+= \g[t] \oplus \Cset K
\label{pa}
\end{equation}
and its ideal
\begin{equation}
\na_+= t \g[t].
\label{na}
\end{equation}
Clearly $\pa_+/\na_+ \cong \g \oplus \Cset K.$
Later we will also need ``the opposite to $\na_+$ subalgebra'' of $\gt$
\begin{equation}
\na_-= t^{-1} \g[t^{-1}].
\label{na-}
\end{equation}

For any subalgebra $\f$ of $\gt$ 
the component of degree $N$ of $U(\f)$  will be denoted by
$U(\f)^N.$ Recall:

\bde{ssmooth} (Kazhdan--Lusztig) For a given $\gt$-module $V$ set
\[
V(N) = \Ann_{U(\na_+)^{-N}} V, \; N \in \Zset_{>0}.
\]
A $\gt$-module $V$ is called strictly smooth if $\cup_N V(N) =V.$
\ede

Note that $V(N) \sub V$ is naturally a $\g$-module for the 
embedding $\g \hra \gt,$ $x \mt x t^0.$

For any $\gt$-module $W$ the strictly smooth part of it
\begin{equation}
W(\infty)= \cup_N W(N)
\label{sspart}
\end{equation}
is a $\gt$-submodule of $W.$ The functor $ W \mt W(\infty)$ 
in the category of $\gt$-modules is left exact.

The commutation relation \eqref{commut} implies that any $\gt$-module in
$\Aff(\C)_\ka$ is strictly smooth. The following equivalent 
characterization of the categories $\Aff(\C)_\ka$ for 
$\ka \notin \Rset_{\geq0}$ 
was obtained in \cite{Y} and in the case $\C=\Fin$ 
$(\ka \notin \Qset_{\geq 0})$
previously in \cite{KL1}.

\bpr{Affcc} Assume that $\C$ is a category of finite length $\g$-modules 
and is closed under tensoring with finite dimensional $\g$-modules.
If $\ka \notin \Rset_{\geq 0}$ then all modules in $\Aff(\C)_\ka$ have 
finite length.
The category $\Aff(\C)_\ka$ coincides with the category
of finitely generated, strictly smooth $\gt$ modules of central charge 
$\ka - h\spcheck$ for which
\[
V(N) \in \C,
\]
considered as a $\g$-module, for all $N \in \Zset_{> 0}.$ 
\epr

For a given $\g$-module $M$ consider the induced $\gt$-modules
\begin{equation}
\label{Ind}
\Ind(M)_\ka = U(\gt) \otimes_{U(\pa_+)} M
\end{equation}
where $\pa_+$ acts on $M$ through the quotient 
$\pa_+ \thra \pa_+/\na_+ \cong \g \oplus \Cset K$ 
and the central element $K$ acts by $\ka - h \spcheck.$ 

If $M$ is an irreducible $\g$-module and $\ka \neq 0$ then $\Ind(M)_\ka$ 
has a unique irreducible quotient, to be denoted by $\Irr(M)_\ka.$
Both modules belong to $\Aff(\C)_\ka$ if $M \in \C$ and the modules of the 
latter type exhaust all irreducible modules in $\Aff(\C)_\ka.$

For central charge $\ka - h\spcheck,$ 
$\ka \notin \Rset_{\geq 0}$
the categories 
$\wh{\H}_\ka,$ $\wh{\O}'_\ka$ have natural 
duality functors, to be denoted by $D(.):$
\begin{equation}
D(V)= (V^d)^\sharp(\infty), \quad V \in \wh{\H}_\ka \; \mbox{or} \; 
\wh{O}'_\ka.
\label{dual-a}
\end{equation}
Given a $\gt$-module $W$ the notation 
$W^\sharp$ stays for its twist by the
automorphism of $\gt$
\[
(x t^k)^\sharp= x (-t)^k, \; x \in \g, \; k \in \Zset; \; 
K^\sharp = - K.
\]

Assume that the $\gt$-module $V \in \wh{\H}_\ka$ or $\wh{\O}'_\ka$ is 
given by \eqref{Sugdecom} with $V^\xi \in \H$ or $\O',$ 
respectively. Denote
\[
(V^d)^\xi = \{ \eta \in V^d \mid \eta|_{V^{\xi'}} = 0
\; \mbox{for} \; \xi' \neq \xi \}.
\]
The underlining space of the dual module $D(V)$ is
the following subspace of
of $V^d$
\begin{equation}
D(V) = 
\bigoplus_{\xi \colon \xi - \xi_i \in \Zset_{\geq 0}}
(V^d)^\xi
\label{Dual}
\end{equation}
and the generalized $L_0$-eigenspaces of $D(V)$ are 
\begin{equation}
\label{dxi}
D(V)^\xi = (V^d)^\xi,
\end{equation}
see \cite{KL1, Y}. The exactness of the functor $D$ on 
$\wh{\H}_\ka$ and 
$\wh{O}'_\ka$ follows from the first equality.

The functors $D$ on $\wh{\H}_\ka$ and $\wh{O}'_\ka$ 
restrict further to a duality functor in $\KL_\ka$ 
which in a more simple way is given by
\begin{equation}
D(V)= (V^\st)^\sharp(\infty), \quad V \in \KL_\ka.
\label{dual-kl} 
\end{equation}
\subsection{Structure of the categories $\wh{O}_\ka$}
Denote by $M_\la$ the highest weight module for $\g$ 
corresponding to highest weight $\la \in \h^\st,$
with respect to 
to the choice of Cartan and Borel subalgebras,
made in Section 2.1.
Denote by $L_\la$ its unique irreducible quotient.
Consider the Cartan and Borel subalgebras
\[ 
\hat{\h}=\h \oplus \Cset K \oplus \Cset d, \;
\wh{\b}_+= \wh{\h} \oplus \n_+ \oplus t \g[t]
\]
of the extended affine Kac--Moody algebra $\hat{\g},$ 
defined in Section 2.3. 
Recall from \cite{Kac} the notation
\begin{align}
&\de \in \wh{\h}^*, \; \de|_{\h}=0, \, \de(K)=0, \, \de(d)=1,  
\nn \\
&\La_0 \in \wh{\h}^*, \; \La_0|_{\h}=0, \, \La_0(K)=1, \, \La_0(d)=0.
\nn
\end{align}  

The induced module $\Ind(M_\la)_\ka,$ viewed as a module over the extended
Kac--Moody algebra by setting $d= -L_0$ is isomorphic to the highest
weight module with highest weight
\begin{equation}
\label{aw}
\la + (\ka - h\spcheck) \La_0 - \frac{(\la, \la+ 2 \rho)}{2 \ka} \de.
\end{equation}

Computing the determinant of the Shapovalov form on a highest weight
module and using the Jantzen filtration, Kac and Kazhdan described the 
structure of Verma modules for
all Kac--Moody algebras by an extension of Jantzen's argument 
\cite{J}. We will state the result for affine 
Kac--Moody algebras in terms of the modules $\Ind(M_\la)_\ka$ rather
than the affine weights \eqref{aw}.

\bde{star} One says \cite{KK} 
that pair $[\la, \mu] \in \h^\st \times \h^\st$ 
satisfies the condition $(\st)$ if there exists a sequence of roots
$\{ \be_i \}_{i=1}^k \sub \De$ and a sequence of nonnegative integers
$\{ m_i \}_{i=1}^k \sub \Zset_{\geq 0}$ $(\be_i \in \De_+$ if $m_i=0)$ 
such that
$\{ \la_i = \la - \be_1 - \ldots - \be_i \}_{i=0}^k$ satisfies
\begin{equation}
\label{st}
\la_{i+1} = r_{\be_i}(\la_i) + \ka m_i \be_i\spcheck \; 
\mbox{and} \;
(\la, \be_i\spcheck) + 2\ka m_i/ |\be_i|^2 \in \Zset_{>0}, 
\; i=1, \ldots, k,
\end{equation}
$\la_k =\mu.$
\ede

The second condition in \eqref{st} means that 
$\la_{i+1} - \la_i \in Q$ which can only happen if 
$\la_{i+1} - \la_i = n_i \be_i$ for some positive integer $n_i.$  

\bth{KK} (Kac--Kazhdan)
The module $\Irr(L_\la)_\ka$ is a subquotient of 
$\Ind(M_\mu)_\ka$ if only if the pair $[\la + \rho , \mu + \rho]$ 
satisfies the condition $(\st).$
\eth

Denote by $\wh{\O}_{\ka, \ig}$ the category of $\gt$-modules,
satisfying the conditions for $\wh{\O}_\ka$ except the condition
to be finitely generated. Each $\gt$-module 
$V \in \wh{\O}_{\ka, \ig}$ has a filtration by $\gt$-submodules
\begin{equation}
\label{Ofilt}
0 = W_0 \sub W_1 \sub \ldots \sub V
\end{equation}
such that each subquotient $W_n/W_{n-1}$ is isomorphic 
to a quotient of $\Ind(L_{\la_n-\rho})_\ka$ for some 
$\la_n \in \h^\st,$ see \cite{DGK}. The modules in $\wh{\O}_\ka$
have finite filtrations with this property.

Consider the equivalence relation on $\h^\st$ induced by $\la \sim \mu$
if the pair $[\la, \mu]$ satisfies the condition $(\st).$ The equivalence 
class of an element $\la \in \h^\st$ will be denoted by $\hchi_\la.$

For a given $\hchi \in (\h/\sim)$ define $\wh{\O}_{\ka, \ig}^{\hchi}$ 
to be 
the subcategory of $\wh{\O}_{\ka, \ig},$ consisting of $\gt$-modules 
having a filtration of the type \eqref{Ofilt} with the property that for 
all subquotients $W_n/W_{n-1},$ $\la_n \in \hchi.$ 

This construction gives rise to the subcategories 
$\wh{\O}^{\hchi}_\ka = \wh{\O}_\ka \cap \wh{\O}^{\hchi}_{\ka, \ig}$
of $\wh{\O}_\ka.$ In the case $\ka \in \Qset_{<0}$ all modules
in $\wh{\O}_\ka$ have finite length and $\wh{\O}^{\hchi}_\ka$ consists
of those modules having a composition series with subquotients
isomorphic to $\Irr(L_\la)_\ka$ for $\la + \rho \in \hchi.$ 

\bth{blockO} (Rocha-Caridi--Wallach \cite{RW} and 
Deodhar--Gabber--Kac \cite{DGK}) 

The category $\wh{\O}_{\ka, \ig}$ has the block 
decomposition
\[
\wh{\O}_{\ka, \ig} = 
\bigoplus_{\hchi \in \h/\sim} \wh{\O}_{\ka, \ig}^{\hchi}.
\]

Moreover
\[ 
\Ext^n_{\wh{\g}, \wh{\h}}(V_1, V_2) =0, \; n
\in \Zset_{\geq 0}
\]
if $V_i \in \wh{\O}^{\hchi_i}_\ka,$ $i=1,2$ and 
$\hchi_1 \neq \hchi_2.$ Here $\Ext^n_{\wh{\g}, \wh{\h}}$
refer to the $\Ext$ groups in the category of $\wh{\g}$-modules
of central charge $\ka -h\spcheck$ which are $\wh{\h}$ locally
finite and semisimple. 
\eth
\sectionnew{Affine Jacquet functors and vanishing of $\Hom$ groups for 
induced affine Harish-Chandra modules}
Throughout this section we will assume
\[
\ka \notin \Rset_{\geq 0},
\]
except in subsection 3.3 where we will generalize most of the 
results to the case of arbitrary noncritical central charge. 
\subsection{Affine Jacquet functors}
Let $V \in \wh{\H}_\ka.$ Define the affine Jacquet functor by
\begin{equation}
\label{ajacquet}
\ja(V) = j(V)^\sharp(\infty).
\end{equation}
Here $j(V)$ refers to the Jacquet module of $V$ considered as a 
$\g$-module under the natural inclusion $\g \hra \gt.$

Since $(\n_0)_\Cset \sub \g$ acts locally nilpotently on $\gt$ the 
subspace $j(V) \sub V^\st$ is a $\gt$-submodule and thus
$j(V)$ is naturally a $\gt$-module of central charge 
$\ka - h\spcheck.$
 
Assume that the module $V$ has the form \eqref{Sugdecom} with 
$V^\xi \in \H.$ The $\g$-submodule
\[
(V^\st)^\xi :=
\{ \eta \in V^\st \mid 
\eta(V^{\xi'})=0 \; \mbox{if} \; \xi' \neq \xi \}
\]
of $V^\st$ is naturally identified with $(V^\xi)^\st.$ 
Under this identification the $\g$-module
\begin{equation}
\label{jj}
\{ \eta \in V^\st \mid  
\eta(V^{\xi'})=0 \; \mbox{if} \; \xi' \neq \xi, \; 
\eta \in \cup_N \Ann_{\n_0^N} V^\st \}
\end{equation}
is naturally identified with $j(V^\xi).$ By abuse of notation we will
denote the $\g$-submodule \eqref{jj} of $V^\st$ by $j(V^\xi).$ Then as a 
$\g$-submodule of $(V^\st)^\sharp(\infty),$ the affine Jacquet module 
\eqref{ajacquet} $\ja(V)$ is
\begin{align}
\label{jjeq}
\ja(V)&=
\bigoplus_{\xi \colon \xi-\xi_i \in \Zset_+} j(V^\xi) 
\\ 
\nn
&\sub \bigoplus_{\xi \colon \xi-\xi_i \in \Zset_+} (V^\xi)^\st
= (V^\st)^\sharp(\infty).
\end{align}
The fact that the rhs of \eqref{jjeq} 
is a subset of the lhs is obvious. The opposite 
inclusion follows from the following lemma proved analogously to 
\cite[Proposition 2.21]{KL1}.

\ble{nn} Let $V \in \wh{\H}_\ka$ then:

(i) There exists a real number $\zeta$ such that the map 
$\g \otimes V^\xi \ra V^{\xi +1}$ given by 
$x \otimes v \mt (x t^{-1})v$
is an isomorphism for $\Re \xi \geq \zeta.$

(ii) For any positive integer $N$ there exists 
$\zeta \in \Rset$ such that
$\oplus_{\Re \xi \geq \zeta} V^\xi \sub U(\na_-)^{-N} V.$
\ele

Eq. \eqref{jjeq} implies that the affine Jacquet functor 
is exact and faithful. It can also be used as a definition of the 
affine Jacquet functor. In this definition of $\ja$  
the (standard) Jacquet functor $j$ is applied only to
the finitely generated, admissible $(\g, \k)$-modules
$V^\xi$ as opposite to  \eqref{ajacquet} where
the Jacquet functor is applied to the $(\g, \k)$-module 
$V \in \GK$ which is neither admissible, nor finitely generated. 

It is straightforward to show that
\begin{equation}
\label{jaxi}
(\ja(V))^\xi = j(V^\xi).
\end{equation}
The fact that $\ja(V)$ has finite length is not obvious.
It will be deduced from the following Proposition.

\bpr{Jaind} If $M \in \H$ then 
\begin{align}
\label{ja1}
\ja(\Ind(M)_\ka) &\cong D( \Ind(j(M)^d)_\ka ),
\\
\label{ja2}
\ja(D(\Ind(M)_\ka)) & \cong \Ind(j(M^d)_\ka)
\end{align}
where in the left and right hand sides $D$ refers to the duality functor 
in the categories $\wh{\H}_\ka$ and $\wh{\O}'_\ka,$ respectively.
\epr 

\prref{Jaind}, the exactness, and the faithfulness of $\ja$ imply that
$\ja(V)$ is a $\gt$-module of finite length for any $V \in \wh{\H}_\ka.$ 
This coupled with \eqref{jaxi} shows that $j : \wh{\H}_\ka \ra 
\wh{O}'_\ka.$ Thus we obtain:

\bpr{Jacquet_a} The affine Jacquet functor $\ja$ given by \eqref{ajacquet} 
or \eqref{jjeq} is faithful and exact (contravariant) functor from 
$\wh{\H}_\ka$ to $\wh{O}'_\ka.$  
\epr

\noindent
{\em{Proof of \prref{Jaind}}} We will show \eqref{ja1}. Eq. 
\eqref{ja2} can be proved analogously. Since the Jacquet functor $j$ and 
the induction functor 
$\Ind(.)_\ka \colon \H \ra \wh{\H}_\ka$ commute with 
finite direct sums, it is sufficient to prove the statement 
of the lemma
for a module $M \in \H^{\chi_\la},$ $\la \in \h^\st.$ Because of 
\eqref{commut}
\begin{equation}
\label{IndM}
\Ind(M)_\ka = \bigoplus_{N \in \Zset_{\geq 0}} 
\Ind(M)_\ka^{\phi(\la) +N}
\end{equation}
where
\begin{equation}
\label{ph}
\phi(\la) = (|\la|^2 - |\rho|^2)/\ka,
\end{equation}
recall \eqref{chiOm}. Under the natural isomorphism of 
$\g$-modules
\[
\Ind(M)_\ka \cong U(\na_-) \otimes_\Cset M,
\]
$\Ind(M)_\ka^{\phi(\la) +N} \sub \Ind(M)_\ka$ is mapped to
$U(\na_-)^{N} \otimes_\Cset M \sub U(\na_-) \otimes_\Cset M.$ 
Here $\g$ acts on $U(\na_-)$ by the adjoint action.
We obtain the chain of natural isomorphisms of $\g$-modules
\begin{align}
\ja(\Ind(M)_\ka) &\cong \bigoplus_{N \in \Zset_{\geq 0}} 
j( U(\na_-)^{N} \otimes_\Cset M) \cong
\bigoplus_{N \in \Zset_{\geq 0}} (U(\na_-)^{N})^\st \otimes_\Cset j(M) 
\cong
\label{isom} \\
&\cong \bigoplus_{N \in \Zset_{\geq 0}} (U(\na_-)^{N} 
\otimes_\Cset j(M)^d)^d 
\cong D(j(M)^d).
\nn
\end{align}

The first isomorphism comes from \eqref{jjeq}, the second is an
application of \prref{finJa}, and the forth one is from \eqref{Dual}.

It is straightforward to check that \eqref{isom} is an isomorphism
of $\gt$-modules.
\qed

\bco{fil} If $M \in \H^{\chi_\la},$ $\la \in \h^\st$ then
$\ja(\Ind(M)_\ka)$ has a filtration by Weyl modules 
$\Ind(L_{-\la' + \rho})_\ka$ for $\la'_i \in W \la.$ 
\eco
\subsection{Structure of induced Harish-Chandra modules}
The main result of this subsection is the following analog of 
Kac--Kazhdan's \thref{KK} for the affine Harish-Chandra categories.

\bth{HomH} If $M_i \in \H,$ $i=1,$ $2$ are two irreducible
$\gt$-modules then $\Irr(M_2)_\ka$ is a subquotient of $\Ind(M_1)_\ka$
only if $M_1,$ $M_2$ have infinitesimal characters $\chi_{\la_1}$ and 
$\chi_{\la_2}$ for a pair $[\la_1,\la_2] \in \h^\st \times \h^\st$ 
satisfying the condition $(\st)$ from \deref{star}.
\eth  
\begin{proof} Assume that $\Irr(M_2)_\ka$ is a subquotient of 
$\Ind(M_1)_\ka.$ Then
\[
\Hom (\Ind(M_2)_\ka, \Ind(M_1)_\ka/V_1) \neq 0 
\]
for some $\gt$-submodule $V_1$ of $\Ind(M_1)_\ka.$
The exactness and the faithfulness of the affine
Jacquet functor imply that 
\begin{equation}
\label{hn}
\Hom (\ja D(\Ind(M_2)_\ka), \ja D(\Ind(M_1)_\ka)/\ja D(V_1)) \neq 0.
\end{equation}
Applying \prref{Jaind}, we see that
$\Irr(j(M_2^d))_\ka$ is a subquotient of $\Ind(j(M_1^d))_\ka.$
Thus $\Irr(L_{\la_2})_\ka$ is a subquotient of $\Ind(L_{\la_1})_\ka$ 
for some $\la_i \in \chi(M_i),$ $i=1,2.$ 
Then \coref{fil} and  
Kac--Kazhdan's \thref{KK} imply the statement of the proposition.
\end{proof}
\subsection{The general case of noncritical central charge}

The affine Jacquet functor \eqref{ajacquet}
and most of its properties can be extended 
to the case $\ka \neq 0.$ This in particular gives
extensions of \thref{HomH} to all noncritical levels
$\ka \neq 0,$ and of \prref{Affcc}  to the case
$\ka \notin \Qset_{\geq 0}.$

In this section we will assume that the central charge
$\ka - h\spcheck$ is noncritical, i.e. that
\[
\ka \neq 0.
\]
Introduce the categories $\wh{\H}_{\ka, \ig}$ of 
$\gt$-modules of central charge $\ka - h\spcheck,$
satisfying the same properties as the modules in $\wh{\H}_\ka$ 
except the condition to be finitely generated.

There are natural duality functors in the categories
$\wh{\H}_{\ka, \ig}$ and $\wh{\O}'_{\ka, \ig}$ (defined in
section 2.4), extending the duality functors $D$ in 
$\wh{\H}_\ka$ and $\wh{\O}'_\ka,$ defined by
\begin{equation}
\label{Dig}
D(V) = ((V^d)^\sharp(\infty) )^{\Cset[L_0]-fin}, \quad
V \in \wh{\H}_{\ka, \ig} \;
\mbox{or} \;  \wh{\O}'_{\ka, \ig}.
\end{equation}
If the $\gt$-module 
$V \in \wh{\H}_{\ka, \ig}$ or
$\wh{\O}'_{\ka, \ig}$ has the decomposition
\eqref{Sugdecom}, then as a submodule of $(V^d)^\sharp$ the dual
module $D(V)$ is given by
\begin{equation}
D(V) = \bigoplus_\xi (V^d)^\xi
\label{Dualig}
\end{equation}
in the notation \eqref{dxi}.

One can define an extension of the affine Jacquet
functor \eqref{ajacquet} to a functor from the category
$\wh{\H}_{\ka, \ig}$ to the category $\wh{\O}'_{\ka, \ig}$
by
\begin{equation}
\label{ajacquetig}
\ja(V) = (j(V)^\sharp(\infty))^{\Cset[L_0]-fin}.
\end{equation}
If the $\gt$-module $V \in \wh{\H}_{\ka, \ig}$ 
has the decomposition \eqref{Sugdecom} then as a submodule of 
$V^\st,$ the affine Jacquet module $\ja(V)$ is given by
\begin{equation}
\label{jjeqig}
\ja(V)=
\bigoplus_{\xi \colon \xi-\xi_i \in \Zset_+} j(V^\xi) \sub
\bigoplus_{\xi \colon \xi-\xi_i \in \Zset_+} (V^\xi)^\st
= ((V^\st)^\sharp(\infty))^{\Cset[L_0]-fin}
\end{equation}
assuming the identification of $j(V^\xi)$ with
the $\g$-submodule \eqref{jj} of $V^\st.$
This implies that
\[
\ja(V) \in \wh{\O}'_{\ka, \ig}, \quad
\mbox{for all} \;  V \in \wh{\H}_{\ka, \ig}
\]
and the validity of an extension of \prref{Jaind} and 
\prref{Jacquet_a}
for the categories $\wh{\H}_{\ka, \ig}.$

\bpr{Jig} (i) The affine Jacquet functor \eqref{ajacquetig}
$\ja \colon \wh{\H}_{\ka, \ig} \ra \wh{\O}'_{\ka, \ig}$
is exact and faithful.

(ii) For any $M \in \H$ we have the isomorphisms 
of $\gt$-modules in the category $\wh{\O}'_{\ka, \ig}$
\begin{align}
\ja(\Ind(M)_\ka) &\cong D( \Ind(j(M)^d)_\ka ),
\nn \\
\ja(D(\Ind(M)_\ka)) &\cong \Ind(j(M^d))_\ka.
\nn
\end{align}
\epr

{}  From part (ii) of \prref{Jig} one obtains the following
extension of \prref{Affcc} and \thref{HomH}.
\bpr{pos} (i) The statement of \thref{HomH} is valid for any
noncritical central charge.

(ii) If $\ka \notin \Qset_{\geq 0}$ then every module in 
the category $\wh{\H}_\ka$ has finite length. In addition the
category $\wh{\H}_\ka$ consists of all finite length, strictly
smooth $\gt$-modules of central charge $\ka - h\spcheck$ for 
which
\[
V(N) \in \H,
\]
considered as $\g$-modules, for all $N \in \Zset_{> 0}.$
\epr
According to \cite[Section 3]{Y}, to prove part (ii) of \prref{pos}
we need to show that $\Ind(M)_\ka$ has finite length 
for any irreducible $\g$-module $M \in \H.$ Because of the faithfulness
of the affine Jacquet functor and part (ii) of \prref{Jig} it is 
sufficient to show that $\Ind(L_{\la-\rho})_\ka$ has finite length
for all $\la \in \h^\st,$ recall the notation from section 2.4. The 
latter follows from \thref{KK} of Kac and Kazhdan. What needs to be 
proved 
is that $\Irr(L_{\mu-\rho})_\ka$ is an irreducible subquotient of
$\Ind(L_{\la-\rho})_\ka$ for finitely many $\mu \in \h^\st.$ 
If the pair $[\la, \mu] \in \h^\st \times \h^\st$ satisfies the condition
$(\st)$ from \deref{star} then 
\begin{equation}
\label{mula}
\mu = \la + \al = w \la  + \ka \be \; 
\mbox{with} \; \al, \be \in Q
\end{equation}
for some element $w \in W.$ Notice that given $w \in W$ there is at
most one weight $\mu \in \h^\st$ that satisfies \eqref{mula} for some
$\al,$ $\be \in Q.$ Indeed if \eqref{mula} is satisfied then
\[
\al + \ka \be = w \la - \la 
\]
which uniquely defines $\al$ and $\be$ since $Q \cap \ka Q =0.$
\sectionnew{Block decomposition of affine Harish-Chandra categories}
\subsection{Relation of $\Ext$ groups for $\gt$ and $\g$-modules}
In this subsection we will assume that the central charge 
$\ka - h\spcheck$ is noncritical, i.e.
\[
\ka \neq 0.
\]

By $\GtK$ and $\PK$ we will denote the categories
of $\Cset$-graded $\gt$ and $\pa_+$-modules, 
respectively, (with respect to the grading \eqref{grading})
of central charge $\ka -h\spcheck$ which are locally finite and semisimple 
when restricted to $\k \sub \gt.$ For each module $V \in \GtK$ or $\PK$
we will denote by $V^\xi$ its degree $\xi(\in \Cset)$ subspace which is 
naturally a $\g$-module in $\GK.$ 

The category $\wh{\H}_\ka$ is canonically embedded in $\GtK$ using the
grading \eqref{Sugdecom} related to the generalized eigenvalue 
decomposition of the Sugawara operator $L_0$ on $V \in \wh{\H}_\ka.$ 
Clearly each morphism in $\wh{\H}_\ka$ preserves this decomposition.

Similarly to the category $\GK,$ the categories $\GtK$ and $\PK$ 
have enough projectives.
Indeed if $V \in \GtK$ is generated as a $U(\gt)$-module
by the $\k$ stable graded subspace $E$ of $V$ then 
\[
U(\gt)/
\all K-\ka + h\spcheck \arr \otimes_{U(\k)} E
\]
is a projective cover of $V.$ Here 
$\all K-\ka + h\spcheck \arr$
denotes the (two-sided) ideal of $U(\gt)$ generated
by $K-\ka + h\spcheck.$ The first factor $U(\gt)$ is 
graded using \eqref{grading}. Analogously one shows 
that $\PK$ has enough projectives. 

By $\Ext_{\gt, \k, \ka}(.,.)$ and 
$\Ext_{\pa_+, \k, \ka}(.,.)$ we will denote the extension
groups in the categories $\GtK$ and $\PK,$ respectively.
For simplicity the $\Hom$ groups in $\GtK$ and $\PK$ will
be denoted by $\Hom_{\gt}$ and $\Hom_{\pa_+}.$

In this subsection we relate $\Ext$-groups between 
induced modules in $\GtK$ to $\Ext$-groups in
$\GK.$

Consider the induction functor
\[
\I(.) \colon \PK \ra \GtK
\]
defined by
\[
\I(P) = U(\gt) \otimes_{U(\pa_+)} P, \quad P \in \PK.
\]
Since it is exact, the canonical isomorphism
\[
\Hom_{\gt} ( \I(P)_\ka, V) \cong 
\Hom_{\pa_+} (P, V), \quad 
P \in \PK, \; V \in \GtK
\]
gives rise to the isomorphisms of the corresponding $\Ext$ groups
\begin{equation}
\label{GtP2}
\Ext_{\gt, \k, \ka} ( \I(P)_\ka, V) \cong 
\Ext_{\pa_+, \k, \ka} (P, V), \quad
P \in \PK, \; V \in \GtK.
\end{equation}
In particular we obtain the isomorphisms
\begin{equation}
\label{GtP}
\Ext_{\gt, \k, \ka}^n ( \Ind(M)_\ka, V) \cong 
\Ext_{\pa_+, \k, \ka}^n (M, V), \quad
M \in \H, \; V \in \GtK
\end{equation}
where in the rhs $M \in \H$ is considered as a $\pa_+$-module in the 
category $\PK$ by setting $\pa_+$ to act on $M$ through the quotient
$\pa_+ \ra \pa_+ /\na_+ \cong \g \oplus \Cset K$ and 
$K$ to act on $M$ by $(\ka -h\spcheck).\id.$ If 
$M = \bigoplus_{\chi_\la \in \h^\st/W} M^{\chi_\la}$ is the 
decomposition of $M$ with respect to \eqref{block} then 
$M^{\chi_\la}\sub M$ is put in degree $\phi(\la),$ recall \eqref{ph}. 
(Recall also the definition \eqref{Ind} of 
the induction functor $\Ind(.)_\ka \colon \GK \ra \GtK.$)

Assume for simplicity that $M \in \H^{\chi_\la}$ for some
$\la \in \h^\st.$
Consider the natural isomorphism 
\begin{equation}
\label{eh}
\Hom_{\pa_+} (M, V) \cong 
\Hom_{\g}(M, (V^{\phi(\la)})^{\na_+}), \quad
M \in \H, \; V \in \GtK,
\end{equation}
where $M$ is thought of as a $\pa_+$ module in $\PK$
as described above.

Given a complex number $\xi,$ denote by $H^q(\na_+,.)^\xi$ the 
$q$-th derived functor of the left exact 
functor $V \mt (V^\xi)^{\na_+}$ from $\PK$ to $\GK,$ for any
$\xi \in \Cset.$ 
A version of the Hochschild--Serre spectral sequence relates
the $\Ext$ groups, corresponding to $\Hom$'s in \eqref{eh}:
\begin{equation}
\label{HS}
E^2_{pq} = \Ext^p_{\g, \k} (M, H^q(\na_+,V)^{ \phi(\la) } ) 
\Rightarrow
\Ext^n_{\pa_+, \k, \ka} (M, V), \quad
M \in \GK, \; V \in \PK.
\end{equation}

Since 
\[
V^{\na_+} = \bigoplus_{\xi \in \Cset} (V^\xi)^{\na_+}, \;
\forall V \in \PK, 
\]
if $H^q(\na_+,.)$ denotes the $q$-th derived 
functor of the functor $V \mt V^{\na_+}$ from $\PK$ to $\GK,$
then
\begin{equation}
\label{Hsum}
H^q(\na_+, .) = \bigoplus_{\xi \in \Cset} H^q(\na_+, .)^{\xi}.
\end{equation}

Combining \eqref{GtP} and \eqref{HS} we obtain:

\bpr{HSss} For all modules 
$M \in \H^{\chi_\la}$ and $V \in \GtK$ one has the spectral sequence
\[
E^2_{pq} = \Ext^p_{\g, \k} (M, H^q(\na_+, V)^{ \phi(\la) }) 
\Rightarrow \Ext^n_{\gt, \k, \ka} (\Ind(M)_\ka, V). 
\]
\epr

{\em{Consider the equivalence relation $\sim$ in $\h^\st$ induced 
by the relation $\la \sim \mu$ if $\mu \in W \la$ or the 
pair $[\la, \mu]$ satisfies the condition $(\st)$ from \deref{star}.}}

\bpr{ExtvanH} If $\la_i \in \h^\st,$ $i=1,2$ and 
$\la_1 \not\sim \la_2,$ then
\begin{equation}
\label{van2}
\Ext^n_{\gt, \k, \ka} ( \Ind(M_1)_\ka, D(\Ind(M^d_2)_\ka)) = 0, 
\end{equation} 
for all $\g$-modules 
$M_1 \in \H^{\chi_{\la_1}}$ and $M_2 \in \H^{\chi_{\la_2}}.$
\epr

We will deduce \prref{ExtvanH} from the following lemma
to be proved at the end of the subsection.

\ble{Hvan} If $M \in \H$ 
\[
H^n(\na_+, D(\Ind(M)_\ka)) =0
\]
for $n >0.$
\ele

\noindent
{\em{Proof of \prref{ExtvanH}.}} \leref{Hvan} shows that the spectral 
sequence from \prref{HSss} degenerates at the $E^2$ term
in the case $V_2 = D(\Ind(M_2)_\ka).$
Therefore
\begin{equation}
\Ext^n_{\gt, \k, \ka} (\Ind(M_1)_\ka, D(\Ind(M_2)_\ka)) \cong
\Ext^n_{\g, \k} (M_1, D(\Ind(M_2)_\ka)^{\na_+}),
\label{2}
\end{equation}
for all $M_1 \in \H^{\chi_{\la_1}}$ $M_2 \in \H^{\chi_{\la_2}}.$ 

\thref{HomH} and part (i) of \prref{pos} imply that any irreducible 
subquotient of $\Ind(M_2^d)_\ka$ is isomorphic to $\Irr(M)_\ka$ for some 
irreducible $\g$-module $M \in \H^{\chi_\la},$ $-\la \sim \la_2.$
Recall that
\begin{equation}
\label{duad}
D(\Ind(M^d)_\ka) \cong \Ind(M)_\ka, \quad \mbox{for all irreducible} 
\; M \in \H
\end{equation}
according to \cite[Proposition 4.1, part 3]{Y}.
Because of this
all irreducible subquotients of $D(\Ind(M_2^d)_\ka$ are of the type
$\Irr(M)_\ka$ for irreducible $\g$-modules $M \in \H^{\chi_\la},$
$\la \sim \la_2.$ Therefore 
\[
(D(\Ind(M_2^d)_\ka))^{\na_+} = D(\Ind(M_2^d)_\ka)(0)
\in \prod_{\la \sim \la_2} \H^{\chi_{\la}}.
\]
Eq. \eqref{2} and \prref{Extvan} imply \eqref{van2}.
\qed

\noindent
{\em{Proof of \leref{Hvan}.}} We follow the strategy of the proof of 
\cite[Proposition 4.7]{DGK} of Deodhar, Gabber, and Kac.  

Because of \eqref{Hsum} we need to show
\begin{equation}
\label{covan}
H^q(\na_+, D(\Ind(M)_\ka))^\xi = 0
\end{equation}
for all $\xi \in \Cset.$ 

Denote by $\Cset$ the trivial $\pa_+$-module in degree 0. We have
the $U(\na_+)$-free resolution of $\Cset[-\xi]$ in $\C_{\pa_+, \k, \ka}:$
\[
\ldots \ra R^1 \ra R^0 \ra \Cset[-\xi] \ra 0
\]
where
\[
R^q = U(\na_+) \otimes_\Cset \wedge^q \na_+ [-\xi]
\]
and $\pa_+$ is defined to act on $R^j$ as follows. The ideal 
$\na_+$ acts by left multiplication on the first term, $\g$ 
acts by adjoint transformation on all factors, and $K$ acts by 
$(\ka-h\spcheck).\id.$ The grading \eqref{grading} is used in the 
definition of the graded components of $R^q.$ By $(.)[-\xi]$ we 
denote the standard operator of shift of the grading 
by $\xi \in \Cset.$

The group $H^q(\na_+, D(\Ind(M)_\ka))^\xi$ is the $q$-th 
homology group
of the associated complex 
$\{ \Hom(R^\bullet, D(\Ind(M)_\ka))) \}.$ 
Note that the latter consists of finite dimensional vector spaces.

One has the isomorphism of vector spaces
\[
\Hom(R^\bullet, D(\Ind(M)_\ka)))
\cong
\Hom_{gr}(\wedge^\bullet \na_+, D(\Ind(M)_\ka)))
\]
where $\Hom_{gr}$ denotes the spaces of degree $0$ homomorphism
of graded vector spaces.

Consider the complex defining the homology groups 
$H_q(\na_-, \Ind(M)_\ka):$ 
\[
0 \ra \Ind(M)_\ka \ra C_0 \ra C_0 \ra \ldots 
\]
where 
\[
C_q = \wedge^q \na_- \otimes_\Cset \Ind(M)_\ka.
\]
It is a complex of $\Cset$ graded vector spaces using the grading 
\eqref{grading} of $\na_-$ and the grading \eqref{Sugdecom} of
$\Ind(M)_\ka.$ The differential preserves this grading. 
Denote the component of degree $\xi \in \Cset$ of $C_q$ by
$C_q^\xi.$ Note that it is a finite dimensional vector space.

There is a natural isomorphism of complexes 
\[
\{ \Hom_{gr}(\wedge^\bullet \na_+, D(\Ind(M)_\ka))) \} 
\cong 
\{ 
\left[ (\wedge^\bullet \na_- \otimes_\Cset \Ind(M)_\ka )^\xi
\right]^\st 
\}
\]
by $\vph \mt \left[ \om \otimes v 
\mt \vph(\om^\sharp)(v) \right],$ $\om \in \wedge^\bullet \na_+,$
$v \in \Ind(M)_\ka.$

Since $H_q(\na_-, \Ind(M)_\ka)=0$ for $q>0$ $(\Ind(M)_\ka$ is a free
$\na_-$-module) we obtain \eqref{covan}.
\qed
\subsection{Block decomposition of $\wh{\H}_\ka$ for
$\ka \not{\in} \Qset_{\geq 0}$}
In this subsection we will assume that
\[
\ka \notin \Qset_{\geq 0}.
\]

For any $\hchi_\la \in (\h^\st/\sim)$ denote by 
$\wh{\H}^{\hchi_\la}_\ka$ the full subcategory of $\wh{\H}_\ka$ consisting 
of $\gt$-modules having composition series with subquotients 
$\Irr(M)_\ka$ for irreducible $\g$-modules with infinitesimal character
$\chi_{\la'}$ for some $\la' \sim \la.$

\thref{HomH} implies
\begin{equation}
\label{indblock}
\Ind(M)_\ka \in \wh{\H}^{\hchi_\la}_\ka, \quad
\mbox{if} \; 
\chi(M)=\chi_\la.
\end{equation}

The main result in this section is the following generalization
of the important \prref{Extvan} on vanishing of $\Ext_{\g, \k}$-groups 
between different blocks in the standard Harish-Chandra category.

\bth{VExt} If $\hchi_1, \hchi_2 \in (\h^\st/\sim)$ are such that
\begin{equation}
\label{lab}
\hchi_1 \neq \hchi_2
\end{equation}
then
\begin{equation}
\label{hmhm}
\Ext^n_{\gt, \k, \ka} (V_1, V_2) =0
\end{equation}
for all $V_1 \in \wh{\H}^{\hchi_1}_\ka$ and  
$V_2 \in \wh{\H}^{\hchi_2}_\ka,$ $n \in\Zset_{\geq 0}.$
\eth

The special case of \thref{VExt} for vanishing of $\Ext^1$ groups 
implies a block decomposition of the categories
$\wh{\H}_\ka$ for $\ka \notin \Qset_{\geq 0}$ which is a generalization of
the block decomposition \thref{blockO} of Deodhar--Gabber--Kac and 
Rocha-Caridi--Wallach for the BGG category $\O$ for Kac--Moody
algebras.

\bpr{block} For $\ka \notin \Qset_{\geq 0}$ the category $\wh{\H}_\ka$
is a direct sum of its subcategories $\wh{\H}_\ka^{\hchi}:$
\[
\wh{\H}_\ka \cong 
\bigoplus_{\hchi \in (\h^\st/\sim)}
\wh{\H}_\ka^{\hchi}.
\]
\epr

We will prove \thref{VExt} in two steps. First we will prove an auxiliary 
lemma.

\ble{part} Assume that $\la_1, \la_2 \in \h^\st$ 
and 
\begin{equation}
\label{la12}
\la_1 \not\sim \la_2
\end{equation}
If $M_1 \in \H^{\chi_{\la_1}}$ 
and $V_2 \in \wh{\H}^{\hchi_{\la_2}}$ then 
\begin{equation}
\label{vv1}
\Ext^n_{\gt, \k, \ka} (\Ind(M_1)_\ka, V_2) =0
\end{equation}
for all $n \in \Zset_{\geq 0}.$
\ele

\begin{proof} We will show \eqref{vv1} by induction 
on $n.$ The case
$n=0$ follows from \thref{HomH}. Assume the validity of \eqref{vv1}
for some $n \in \Zset_{\geq 0}.$ Every $\gt$-module in 
$\wh{\H}^{\hchi_{\la_2}}$ has a finite composition series 
with subquotients of the type $\Irr(M)_\ka$ for irreducible
$\g$-modules $M \in \H$ with infinitesimal character
$\chi(M)=\chi_{\la'}$ such that $\la' \sim\la_2.$
Using the long exact sequence for the $\Ext$ groups one see 
that to prove \eqref{vv1} for $n+1$ one only needs to show
that
\[
\Ext^n_{\gt, \k, \ka} (\Ind(M_1)_\ka, \Irr(M_2)_\ka) =0
\]
for $M_1 \in \H^{\chi_{\la_1}}$ and $M_2 \in \H^{\chi_{\la_2}},$
assuming \eqref{la12}.

Denote the maximal, nontrivial $\gt$-submodule of $\Ind(M_2)_\ka$
by $X(M_2^d)_\ka.$ Then one has the exact sequence of $\gt$-modules
\[
0 \ra X(M_2^d)_\ka \ra \Ind(M_2^d)_\ka \ra \Irr(M_2^d)_\ka \ra 0.
\]
Applying the exact duality functor $D$ and using the fact
that $D(\Irr(M_2^d)_\ka) \cong \Irr(M_2)_\ka,$ see 
\eqref{duad}, we get the exact sequence
\begin{equation}
\label{ex}
0 \ra \Irr(M_2)_\ka \ra D(\Ind(M_2^d)_\ka) \ra 
D(X(M_2^d)_\ka) \ra 0.
\end{equation}
This finally gives the exact sequence of $\Ext$ groups
\begin{align}
\Ext^n_{\gt, \k, \ka} (\Ind(M_1)_\ka, D(X(M_2^d)_\ka) ) &\ra
\Ext^{n+1}_{\gt, \k, \ka} (\Ind(M_1)_\ka, \Irr(M_2)_\ka ) \ra
\\
\nn
&\ra \Ext^{n+1}_{\gt, \k, \ka} (\Ind(M_1)_\ka, D(\Ind(M_2^d)_\ka) ).
\nn
\end{align}
The induction implies that the first group vanishes because 
$D(X(M_2^d)_\ka)$ is a quotient of 
$D(\Ind(M_2^d)_\ka) \in \wh{\H}^{\hchi_{\la_2}}.$
The third group vanishes, as proved in \prref{ExtvanH}.
This implies that the middle group vanishes and completes the 
induction argument.
\end{proof} 

\noindent
{\em{Proof of \thref{VExt}.}} Again we use induction on $n.$ The case
$n=0$ is trivial. Assume the validity of the statement in \thref{VExt}
for an integer $n.$ Since every module in the categories $\wh{\H}_\ka$
has finite length for $\ka \notin \Qset_{\geq 0},$ to show that 
\eqref{hmhm} holds for $n+1,$ we only need to prove 
\begin{equation}
\label{no}
\Ext^n_{\gt, \k, \ka} ( \Irr(M_1)_\ka, V_2) = 0
\end{equation}
for an irreducible $\g$-module $M_1 \in \H$ with infinitesimal character
$\chi(M_1) = \chi_{\la_1} \in \h^\st/W$ such that 
$\hchi_{\la_1} \neq \hchi_2.$
Using the exact sequence 
\[
0 \ra X(M_1)_\ka \ra \Ind(M_1)_\ka \ra \Irr(M_1)_\ka \ra 0
\]
we obtain the exact sequence
\begin{align}
\Ext^n_{\gt, \k, \ka} (X(M_1)_\ka, V_2) &\ra
\Ext^{n+1}_{\gt, \k, \ka} (\Irr(M_1)_\ka, V_2) \ra
\nn
\\
&\ra \Ext^{n+1}_{\gt, \k, \ka} (\Ind(M_1)_\ka, V_2) 
\nn
\end{align}
The first term vanishes due to the induction hypothesis. 
The third term vanishes because of \leref{part}. This implies
\eqref{no} and completes the proof of \thref{VExt}.
\qed
\subsection{A block decomposition of $\wh{\H}_{\ka, \ig}$
for arbitrary noncritical central charge}
In this section we extend the block decomposition of the 
categories $\wh{\H}_\ka$ from section 4.2 
to the case of arbitrary noncritical central charge.
Even more generally we provide such a decomposition for 
the categories $\wh{\H}_\ka$ of not necessarily finitely 
generated affine Harish-Chandra modules. Note that they are 
naturally embedded in the categories $\GtK,$ defined in section 4.1. 

It is easy to show by induction that:

{\em{ $(\bullet)$ Every module $V$ in $\wh{\H}_{\ka, \ig}$ has an 
increasing 
filtration 
\[
0= W_0 \sub W_1 \sub \ldots \sub V
\]
for which the subquotients $W_i/W_{i-1}$ are isomorphic to quotients of 
induced modules $\Ind(M_i)_\ka$ for some irreducible $\g$-modules 
$M_i \in \H.$}} 
 
The subcategory $\wh{\H}_\ka$ of $\wh{\H}_{\ka, \ig}$ consists
of those $\gt$-modules in $\wh{\H}_{\ka, \ig}$ that have
a finite filtration with the property $(\st).$ If $\ka \in \Qset_{>0}$
the modules in $\wh{\H}_\ka$ have in general infinite length.

Given any $\hchi_\la \in (\h^\st/\sim)$ for some $\la \in \h^\st$
denote by $\wh{\H}^{\hchi_\la}_{\ka, \ig}$ the subcategory of 
$\wh{\H}_{\ka, \ig}$ consisting of $\gt$-modules possessing a filtration
with the property $(\bullet)$ such that
\begin{equation}
\label{M}
\chi(M_i)=\chi_{\la'} \;
{\mbox{for some}} \; \la' \sim \la.
\end{equation}
\thref{HomH} and part (i) of \prref{pos}
imply that $\wh{\H}^{\hchi_\la}_{\ka, \ig}$ can be characterized as
the subcategory of $\wh{\H}_\ka$ consisting of $\gt$-module for which 
all irreducible subquotients are isomorphic to $\Irr(M)_\ka$ for 
some irreducible $\g$-module $M \in \H$ satisfying \eqref{M}. 

The subcategories $\wh{\H}^{\hchi_\la}_\ka = \wh{\H}_\ka 
\cap \wh{\H}^{\hchi_\la}_{\ka, \ig}$ of $\wh{\H}_\ka$ consist
of those $\gt$-modules $V \in \wh{\H}_\ka,$ having a finite
filtration by $\gt$-modules with subquotients that are isomorphic
to quotients of induced modules $\Irr(M)_\ka$ for irreducible
$\g$-modules $M \in \H$ with infinitesimal character as in \eqref{M}.  
In the case $\ka \notin \Qset_{\geq 0}$ the definition of 
$\wh{\H}^{\hchi_\la}_\ka$ is consistent with the one from 
section 4.2.

The main result of this section is the following generalization
of the basic decomposition theorem of Deodhar--Gabber--Kac 
and Rocha-Caridi--Wallach for the BGG category $\O$ for Kac-Moody 
algebras.

\bth{blockr} For any noncritical central charge
the category $\wh{\H}_{\ka, \ig}$ is a direct 
sum of its subcategories $\wh{\H}^{\hchi}_{\ka, \ig}:$
\[
\wh{\H}_{\ka, \ig} = \bigoplus_{\hchi \in (\h^\st/\sim)}
\wh{\H}^{\hchi}_{\ka, \ig}.
\]
\eth

\thref{blockr} follows from \prref{fV} below as 
\cite[Theorem 4.2]{DGK} of Deodhar, Gabber, and Kac.

First we state two corollaries of \thref{blockr}. There exists a coarser 
block decomposition of the affine Harish-Chandra categories $\wh{\H}_\ka$ 
which more closely resembles the finite dimensional case. To state it, 
consider the equivalence relation on $\h^\st,$ defined by 
\[
\la \approx \mu \; \mbox{if} \; \mu \in W \la + Q\spcheck \; \mbox{and} \;
\mu - \la \in Q.
\]
\bco{block2} For any noncritical central charge the affine Harish-Chandra 
category has the block decomposition
\[
\wh{\H}_\ka = \bigoplus_{\hchi \in (\h^\st/\approx)} 
\wh{\H}_\ka^{\hchi}
\]
where for $\hchi \in (\h^\st/\approx),$ $\wh{\H}_\ka^{\hchi}$
is the subcategory of $\wh{\H}_\ka$ consisting of affine Harish-Chandra
modules possessing a finite filtration with subquotients isomorphic to 
quotients of induced modules $\Ind(M)_\ka$ for $\chi(M) = \chi_\la,$ 
$\la \in \hchi.$
\eco

This block decomposition is even more explicit in the case of
rational central charge. 
Define $\wh{\H}^{\int}_\ka$ to be the subcategory of $\wh{\H}_\ka$ 
consisting of $\gt$-modules $V$ for which all 
$\g$-submodules $V^\xi$ have composition series with subquotients
with infinitesimal characters in $P/W$ where $P$ is the weight lattice of 
$\g.$ In other words $\wh{\H}_\ka^{\int}$ is the subcategory of 
$\wh{\H}_\ka$ consisting of $\gt$-modules possessing a finite filtration 
with subquotients isomorphic to quotients of induced modules 
$\Ind(M)_\ka$ for $\chi(M) = \chi_\la,$ $\la \in P.$ 

\bco{block3} Assume that $\ka = p/q$ for two nonzero, relatively prime 
integers $p$ and $q.$ Then the subcategory $\wh{\H}^{\int}_\ka$ of
$\wh{\H}_\ka$ has the block decomposition
\[
\wh{\H}^{\int}_\ka = 
\bigoplus_{\hchi_\la \in \h^\st/(W \ltimes p Q\spcheck)} 
\wh{\H}_\ka^{\int, \hchi_\la}
\]
where $\wh{\H}^{\int, \hchi_\la}$ is the subcategory of 
$\wh{\H}^{\int}_\ka$ 
consisting of modules possessing a finite filtration with subquotients 
isomorphic to quotients of induced modules $\Ind(M)_\ka$ for 
$\chi(M) = \chi_\mu,$ $\mu \in W\la + p Q\spcheck$.  
\eco

\bre{block3} \coref{block2} and \coref{block3} are easily extended to 
the infinitely generated version (the category 
$\wh{\H}^{\int}_{\ka, \ig}$) of the category 
$\wh{\H}^{\int}_\ka.$ The category 
$\wh{\H}^{\int}_{\ka, \ig}$ is the full subcategory
of $\wh{\H}_{\ka, \ig},$ consisting of $\gt$-modules $V$ for which all 
$\g$-submodules $V^\xi$ have composition series with subquotients
with infinitesimal characters in $P/W.$ We state them only in the 
case of $\wh{\H}^{\int}_\ka$ to avoid extra technical details.
\ere

Now we return to the proof of \thref{blockr}.

\bpr{fV} If $\hchi_1,$ $\hchi_2 \in (\h^\st/\sim)$ 
and $\hchi_1 \neq \hchi_2,$ then 
\[
\Ext^1_{\gt, \k, \ka}(V_1, V_2) = 0, \;
\mbox{for all} \; V_1 \in \wh{\H}^{\hchi_1}_\ka, \;  
V_2 \in \wh{\H}^{\hchi_2}_\ka.
\]  
\epr

For the proof of \prref{fV} we will need the following lemma.

\ble{Wfilt} Let $V \in \wh{\H}_{\ka, \ig}.$ For any real number 
$\zeta$ the module $V$ has a finite filtration by $\gt$-submodules 
\[
0= W_0 \sub W_1 \sub \ldots \sub W_N = V
\]
such that there exists a subset $J \sub \{1, \ldots N\}$
with the following properties.

(i) If $j \in J,$ then
\[
V_j / V_{j-1} \cong \Irr(M_j)_\ka, \quad 
{\mbox{for some irreducible $\g$-modules}} \; M \in \H.
\]

(ii) If $ \notin J,$ then $(V_j/V_{j-1})^\xi=0$ for $\Re \xi \leq \zeta.$
\ele

For a $\gt$-module $V \in \wh{\H}_{\ka, \ig}$ and a real number 
$\zeta$ define the $\g$-module
\begin{equation}
\label{zeta}
V^{\all \zeta \arr} = 
\bigoplus_{\xi \in \Cset \colon \Re \xi \leq \zeta} V^\xi \in \H.
\end{equation}
Note that the above sum is finite. 

We prove {\em{\leref{Wfilt}}} by induction on the length of the 
$\g$-module $V^{\all \zeta \arr}.$

If $V^{\all \zeta \arr}$ is trivial, then the statement of the lemma is 
obvious. Otherwise choose $\xi_0 \in \Cset$ such 
that $\Re \xi_0 \leq \zeta,$
\begin{equation}
\label{zeta2}
V^{\xi_0} \neq 0 \quad \mbox{and} \quad 
V^{\xi_0-n} = 0 \; \mbox{if} \; n \in \Zset_{>0}.
\end{equation}
Let $M$ be an irreducible $\g$-submodule of $V^{\xi_0}.$ Then
$M$ is annihilated by $\na_+$ and we get a canonical homomorphism
$f \colon \Ind(M)_\ka \ra V.$ Denote the image of $f$ by $U_2$
and the image of the maximal, nontrivial submodule of $\Ind(M)_\ka$
under $f$ by $U_1.$ Then 
\[
U_2/U_1 \cong \Irr(M)_\ka
\]
and the lengths of the $\g$-modules $U_1^{\all \zeta \arr}$ and 
$(V/U_2)^{\all \zeta \arr}$ are strictly less than the length 
$V^{\all \zeta \arr}.$ 
This completes the induction argument.  
\qed

We need an auxiliary lemma for the proof of \prref{fV}. 

\ble{aux} Let $M \in \H$ be an irreducible $\g$-module
with infinitesimal character $\chi_\la,$ $\la \in \h^\st.$

(i) If $V_2 \in \wh{\H}_\ka$ is such that 
\[
V^{\all \Re \phi(\la) \arr}_2 = 0
\]
then $\Ext_{\gt, \k, \ka}( \Irr(M)_\ka, V_2) =0,$
recall \eqref{ph}.

(ii) If $V_2 \in \wh{\H}^{\hchi_2}_\ka,$ and $\hchi_\la \neq \hchi_2$ 
then
\[
\Ext^1_{\gt, \k, \ka}(\Ind(M)_\ka, V_2) = 0.
\]
\ele
\begin{proof} Part (i): We need to show that any extension
in $\GtK$ of $\Ind(M)_\ka$ by $V_2$ splits. Consider one such extension
\begin{equation}
\label{se}
0 \ra V_2 \ra W \ra \Ind(M)_\ka \ra 0.
\end{equation}
The functor $V \mt V^{\xi}$ is an exact functor from 
$\GtK$ to $\H$ for any complex number $\xi.$ 
This, combined with the hypothesis, implies that
\[ 
W^{\phi(\la)} \cong \Ind(M)_\ka^{\phi(\la)}
\cong \Ind(M)_\ka (0) \cong M.
\]
Moreover $W^{\phi(\la)-n}=0$ for $n \in \Zset_{>0}$ because of the 
assumption and \eqref{IndM}. Thus $W^{\phi(\la)} \sub W(0)$ and there 
exists a canonical homomorphism from $\Ind(M)_\ka$ to $W$ which defines a 
splitting of \eqref{se}.

Using the long exact sequence for the $\Ext$ groups and \leref{Wfilt},
one sees that to prove {\em{part (ii)}} it suffices to to show that 
\[
\Ext_{\gt, \k, \ka}(\Ind(M)_\ka, \Irr(M_2)_\ka) =0
\]
for any irreducible $\g$-module $M_2 \in \H$ with infinitesimal
character $\chi(M_2)= \chi_{\la_2}$ such that
\[
\la_2 \notin W \la + p Q.
\]
This is done analogously to \leref{part}.
Recall from \eqref{ex} the exact sequence
\[ 
0 \ra \Irr(M_2)_\ka \ra D(\Ind(M^d_2)_\ka) \ra D(X(M^d_2)_\ka) \ra 0.
\]
It gives rise to the exact sequence
\begin{align}
\nn
\Hom(\Ind(M)_\ka, D(X(M^d_2)_\ka)) &\ra 
\Ext^1_{\gt, \k, \ka}(\Ind(M)_\ka, \Irr(M_2)_\ka) \ra
\\
\nn
&\ra \Ext^1_{\gt, \k, \ka}(\Ind(M)_\ka, D(\Ind(M^d_2)_\ka) ).
\end{align}
The first term vanishes because 
$D(X(M^d_2)_\ka) \in \wh{\H}_{\ka, \ig}^{\hchi_{\la_2}}.$ The third term 
vanishes as a consequence of \prref{ExtvanH}. This completes the proof of 
part (ii). 
\end{proof}

\noindent
{\em{Proof of \prref{fV}.}} We will show that
\begin{equation}
\label{ax}
\Ext^1_{\gt, \k, \ka} (U, V_2) =0
\end{equation}
for any $\gt$-modules that is a quotient of an induced module 
$\Ind(M)_\ka$ from an irreducible $\g$-module $M \in \H^{\chi_\la}$
such that $\la$ is in the class of $\hchi_1.$ This implies the 
general case using the long exact sequence for the 
$\Ext$ groups and the fact that any module 
$V_1 \in \wh{\H}_\ka^{\hchi_1}$ has a finite filtration by 
such $\gt$-modules $U.$

Let 
\[
0 \ra X \ra \Ind(M)_\ka \ra U \ra 0
\]
be the long exact sequence defining $U.$ From it we deduce
the exact sequence 
\begin{align}
\nn
\Hom(X, V_2) &\ra 
\Ext^1_{\gt, \k, \ka}(U, V_2) \ra
\\
&\ra \Ext^1_{\gt, \k, \ka}(\Ind(M)_\ka, V_2).
\nn
\end{align}
Since $X \in \wh{\H}^{\hchi_1}$ and $\hchi_1 \neq \hchi_2$ 
the $\Hom$ group vanishes. The last $\Ext$ group vanishes 
due to \leref{aux}. This implies that the middle term vanishes
which proves \eqref{ax}.
\qed
\sectionnew{Compatibility of the affine Jacquet functor with the 
Kazhdan--Lusztig tensor product}
In this section we prove a compatibility relation between the
affine Jacquet functor and the Kazhdan--Lusztig tensor product
which can be considered as an affine (fusion) analog of 
\prref{finJa}. 

Throughout the section we will assume
\[
\ka \notin \Qset_{\geq 0}
\]
and will use the notation from the appendix.
\subsection{Main result}
\bth{JKL} For any two $\gt$-modules $U \in \KL_\ka$ and
$V \in \wh{\H}_\ka$ there exists a canonical isomorphism
\begin{equation}
\label{JKLisom}
\ja (U \dot{\otimes} V) \cong 
D \left[ U \dot{\otimes} D \ja(V)  \right].
\end{equation}
\eth

\thref{JKL} follows from \prref{gtequ} and \prref{Ga2} below.

We will need some notation. For $M \in \H$ set
\begin{equation}
\label{vph}
\vph(M) = \left( M^d \right)^\st.
\end{equation}
It defines an exact functor from the category of Harish-Chandra
modules to the category of $\g$-modules. For $V \in \wh{\H}_\ka$ set
\begin{equation}
\label{Phi}
\Phi(V) = \left[ D(V) \right]^{\st \sharp} (\infty).
\end{equation}
If $V$ has the decomposition \eqref{Sugdecom} with respect to the action
of $L_0,$ then the underlining $\g$-module of $\Phi(V)$ is canonically 
identified with
\[
\bigoplus_{\xi \colon \xi -\xi_i \in \Zset_{\geq 0} }
\vph(V^\xi).
\] 
Clearly $\Phi$ is an exact functor from $\wh{\H}_\ka$ to the 
category of strictly smooth $\gt$-modules of central charge 
$\ka - h\spcheck.$

Given $V \in \H,$ the embedding of $\g$-modules $V(N) \hra V$ induces the 
embedding
\begin{equation}
\label{Vem}
\vph(V(N)) \hra \Phi(V).
\end{equation}
It is easy to see:
\ble{Viso} For any $V \in \wh{\H}_\ka$ the embedding \eqref{Vem} 
is an isomorphism of $\g$-modules.
\ele

For a $\g$-module $M$ define also its $\g$-submodule
\[
j^0(M) = \cup_{N \in \Zset_{\geq 0}} \Ann_{\n_0^N} M.
\]
The Jacquet module of $M \in \H$ is given by $j(M) \cong j^0(M^\st).$
If $V$ is a $\gt$-module then $j^0(V)$ is a $\gt$-module as well.

Fix two $\gt$-modules $U$ and $V$ as in \thref{JKL}.
Set for shortness 
\begin{equation}
\label{WW}
W = U \otimes V.
\end{equation}

The natural embedding $\g \hra \Ga$ by constant functions
induces structures of $\g$-modules on the spaces
$W$ and $W/G_k W$ for all $k \in \Zset_{>0},$ as well as
on their full and restricted duals. Note the canonical isomorphism 
of $\g$-modules
\begin{equation}
\label{cang}
\Ann_{G_k} W^d = (W/G_kW)^d, 
\quad\Ann_{G_k} W^\st = (W/G_kW)^\st.
\end{equation}

\leref{tensorle} and the fact that 
$\H$ is stable under tensoring with finite dimensional 
$\g$-modules imply 
\begin{equation}
\label{WH}
W / G_k W \in \H, \; k \in \Zset_{>0}.
\end{equation}

\ble{ten} Under the above assumptions on the $\gt$-modules 
$U,$ $V,$ and $W$
there exists a canonical isomorphism of $\g$-modules
\[
\left[ \left( \Ann_{G_k} W^d \right)^d \right]^\st
\cong \Ann_{G_k} W^\st,
\]
for all $k \in \Zset_{>0}.$
\ele
\begin{proof} Because of \eqref{WH} we have the canonical 
isomorphism of $\g$-modules
\[
\left( (W/G_k W)^d \right)^d \cong W/G_kW.
\]
Combining this with \eqref{cang}
gives the need canonical isomorphism of $\g$-modules
\[
\left[ \left( \Ann_{G_k} W^d \right)^d \right]^\st \cong
\left[ \left( (W/G_k W)^d \right)^d \right]^\st \cong
(W/G_kW)^\st \cong \Ann_{G_k} W^\st.
\]
\end{proof}

One can define a structure of a strictly smooth $\gt$-module on the
space
\[
\cup_{k \in \Zset_{>0}} \Ann_{G_k} W^\st
\]
in exactly the same way as the $\g$-action on 
\[
T'(U, V) = \cup_{k \in \Zset_{>0}} \Ann_{G_k} W^d
\]
from the appendix.

\bpr{gtequ} In the above notation for $U,$ $V,$ and $W$ 
there exists a canonical isomorphism of $\gt$-modules
\[
\Phi \left[ \cup_{k \in \Zset_{>0}} \Ann_{G_k} W^d \right] \cong
\cup_{k \in \Zset_{>0}} \Ann_{G_k} W^\st.
\]
Under it the $\gt$-submodule $\ja(U \dot{\otimes} V)$ of the first 
module is sent to the $\gt$-submodule 
$j_0(\cup_{k \in \Zset_{>0}} \Ann_{G_k} W^\st)$ of the second one.
\epr

\begin{proof} \leref{Viso}, \leref{ten}, and \eqref{NN} give rise to the 
canonical isomorphisms of $\g$-modules
\begin{equation}
\label{isoo}
\Phi \left[ \cup_{k \in \Zset_{>0}} \Ann_{G_k} W^d \right] (N) \cong
\vph(\Ann_{G_N} W^d) \cong \Ann_{G_N} W^\st.
\end{equation}
One checks directly that they are consistent for different 
$N \in \Zset_{>0}$ and when put together, define an isomorphism
of $\gt$-modules as stated in \eqref{isoo}.

The second part is straightforward.
\end{proof}

Similarly to the above discussion, there is a canonical embedding
of $\g$-modules
\[
\left( D\ja(V) \right)^d  \hra j^0 V^\st.
\]
It gives rise to the embedding of $\Ga$-modules
\[
\left( U \otimes D\ja(V) \right)^d 
\hra U^\st \otimes \left( D \ja(V) \right)^d 
\hra U^\st \otimes j^0 V^\st.
\]
One further has the embedding of $\Ga$-modules
\[
U^\st \otimes j^0 V^\st \hra (U \otimes j^0 V)^\st.
\]

\bpr{Ga2}The above embeddings of $\Ga$-modules
\[
\left( U \otimes D\ja(V) \right)^d \hra
U^\st \otimes j^0 V^\st \hra (U \otimes j^0 V)^\st
\]
gives rise to the isomorphisms of $\Ga$-modules
\begin{align}
\label{ggco}
\cup_{k \in \Zset{>0}} \Ann_{G_k} \left( U \otimes (D\ja(V)) \right)^d 
&\cong
\cup_{k \in \Zset{>0}} \Ann_{G_k} \left( U^\st \otimes j^0 (V^\st)) 
\right) \cong
\\ 
&\cong
\cup_{k \in \Zset{>0}} j_0 \Ann_{G_k} (U \otimes V )^\st. 
\nn
\end{align}
Each of the three spaces in \eqref{ggco} is equipped with a structure
of $\g$-module analogously to $T'(U, V).$
\epr
\subsection{Finiteness of the fusion tensoring of $\wh{\H}_\ka$
by $\KL_\ka$}
Here we modify \thref{JKL} to give another proof of the case of 
Theorem 1.6 from our work \cite{Y} when the reductive subalgebra $\f$ of 
$\g$ is equal to $\k$ (see \thref{Y} below).

For any positive integer $N$ denote by $\O^{'N}$ the truncation 
of the category $\O'$ for $\g,$ consisting of finitely 
generated $\g$-modules $M$ with (generalized) weight space decompositions
\[
M = \bigoplus_{\la \in \h^\st} M^\la 
\]
such that
\[
\prod_{i=1}^N (h_i - \la(h_i))^N m =0, \; \mbox{for all} 
\; m \in M, \; h_i \in \h.
\]
Similarly the categories $\wh{\O}'_\ka$ posses the truncations
$\wh{\O}^{'N}_\ka$ consisting of finitely generated, strictly smooth  
$\gt$-modules $V$ for which
\[
V(N) \in \O^{'N}, \; \forall N \in \Zset_{>0}.
\]
The union of the categories $\wh{\O}^{'N}_\ka$ is the category 
$\wh{\O}'_\ka.$ 

The categories $\wh{\O}^{'N}_\ka$ can be also described in terms 
of the Borel and Cartan subalgebras of the extended affine
Kac--Moody algebra $\wh{\g},$ see section 2.4. Consider the categories of 
finitely generated modules $V$ of central charge $\ka - h\spcheck$ for 
$\wh{\g}$ with the properties:

(i) The subalgebra $\n_+ \oplus t \g[t]$ of $\wh{\g}$ acts 
locally nilpotently on $V.$

(ii) The module $V$ has a generalized weight space decomposition with 
respect to $\wh{\h}$ 
\[
V = \bigoplus_{\wh{\la} \in {\wh{\h}}^\st}
V^{\wh{\la}}, \; \mbox{such that for all} \;
v \in V^{\wh{\la}} \;
\prod_{i=1}^N (h_i - \wh{\la}(h_i))^N v = 0 \;
\forall h_i \in \wh{\h}.
\]

The categories $\wh{\O}^{'N}_\ka$ are essentially the above 
truncations of the generalized BGG category for $\wh{\g}:$
for each module from $\wh{\O}^{'N},$ consider infinitely many copies of
it on which the operator $d$ acts by $const - L_0$ for an arbitrary 
choice of the constant.

\ble{ON} Any strictly smooth $\gt$-module $V$ of central charge $N$
such that 
\begin{equation}
\label{x}
V(k) \in \O^{'N}, \; \forall \; k \in \Zset_{>0}
\end{equation}
has a finite filtration by $\gt$ submodules 
\[
0= W_0 \sub W_1 \sub \ldots \sub W_L=V
\]
with strictly smooth $\gt$-submodules $W_i$ such that
\begin{equation}
\label{Wi}
(W_{i+1}/W_i)(k) \in \O, \; \forall \; k \in \Zset_{>0}. 
\end{equation}
\ele
\begin{proof} For $j\in \Zset_{>0}$ define the subspace $W_j$ of 
$V$ by
\[
V_j = \Span \{ v \in V \mid \; \mbox{there exists} \; \la \in \h^\st
\; \mbox{such that} \; 
\prod_{i=1}^j(h_i - \la(h_i))^i v = 0, \, \forall h_i \in \h \}.
\]
It is a $\gt$-submodule of $V$ because the adjoint action of 
$\h$ on $\gt$ 
is semisimple. For some sufficiently large integer 
$L \leq N . \dim \h$
we have $W_L =V.$

Define the generalized weight spaces of $V$ with respect to $\h$ by
\[
V^\la = \{ v \in V \mid (h - \la(h))^N v = 0 \; 
\mbox{for all} \; h \in \h \},
\]
for $\la \in \h^\st.$
The assumption \eqref{x} implies that 
$V \bigoplus_{\la \in \h^\st} V^\la.$ 

Observe that for each $h \in \h$ the linear transformation
\[
\psi_h(v) = (h - \la(h))v, \; \mbox{if} \; v \in V^\la
\]
is a $\gt$-endomorphism of $V.$

We will show that the $\gt$-module $V/W_1$ satisfies
\[
(V/W_1)(k) \in \O^{'(L-1)}, \; \forall k \in \Zset_{>0}.
\] 
{}   From this \eqref{Wi} follows by induction. 
Moreover it is obvious that
$(V/W_1)(k)$ is spanned by vectors $v$ such that
\[
\prod_{i=1}^{L-1} (h_i - \la(h)) v = 0, \; 
\forall h_i \in \h
\]
and $\n_+$ acts locally nilpotently on $(V/W_1).$ 
So we only need to show that
$(V/W_1)(k)$ are finitely generated $\g$-modules.

Denote the projection $p \colon V \ra V/W_1.$ Then 
\begin{equation}
\label{p}
(V/W_1)(k) \cong p^{-1} \left[(V/W_1)(k) \right] /
\left( p^{-1} \left[(V/W_1)(k) \right] \cap W_1 \right).
\end{equation}
Besides this 
\begin{align}
p^{-1} \left[(V/W_1)(k) \right] &= \{ v \in V \mid
\psi_h(u v)=0, \; \forall u \in U(\na_+)^{-k}, \; h \in \h \}
\label{z}
\\
&=\cap_{h \in \h} \psi^{-1}_h \left[ V(N) \right].
\nn
\end{align}
Fix a basis $h_1,$ $\ldots,$ $h_r,$ of $\h.$ Then \eqref{p} and \eqref{z} 
imply that $(V/W_1)(k)$ is isomorphic to the
image of the $\g$-homomorphism
\[
(\psi_{h_1}, \ldots, \psi_{h_r}) \colon 
p^{-1} \left[(V/W_1)(k) \right] \ra V(k) \oplus \ldots \oplus V(k),
\; (r \; \mbox{times}).
\]
Since $V(k)$ is finitely generated this gives that $(V/W_1)(k)$
is a finitely generated $\g$-module.
\end{proof}

The main ingredient of Kazhdan--Lusztig's finiteness result \cite{KL1}
that $\KL_\ka$ is closed under the fusion tensor product is a 
hard characterization of the category $\wh{\O}_\ka$ only in terms of 
the modules $V(k).$ It is based on the BGG reciprocity for Kac--Moody 
algebras of Rocha-Caridi and Wallach \cite{RW} or equivalently
Soergel's generalized BGG reciprocity \cite{BGS}. Later this 
characterization was extended by Finkelberg to the category 
$\wh{\O}_\ka,$ \cite{F}. 

\bth{KLF} (Kazhdan--Lusztig \cite{KL1}, Finkelberg \cite{F}) If $V$ is a 
strictly smooth $\gt$-module of central charge $\ka - h\spcheck$ such that
\[
V(k) \in \O, \; \forall k \in \Zset_{>0}
\]
then $V$ has finite length and thus $V \in \wh{\O}_\ka.$
\eth

\thref{KLF}, combined with \leref{ON}, implies:

\bpr{StrOO} If $V$ is a strictly smooth $\gt$-module of central charge 
$\ka - h\spcheck$ such that
\[
V(k) \in \O^{'N}, \; \forall k \in \Zset_{>0}
\]
for some fixed integer $N,$ then $V$ has finite length and
$V \in \wh{\O}^{'N}_\ka.$  
\epr

Finally using affine Jacquet functors and \prref{StrOO} we obtain a 
second proof of the following special case of our finiteness
result \cite[Theorem 1.6]{Y} for the Kazhdan--Lusztig tensor product.

\bth{Fin} Assume that $\ka \notin \Qset_{\geq 0}.$ Then the categories
$\wh{\O}^{'N}_\ka,$ $\wh{\O}'_\ka,$ and $\wh{\H}_\ka$ are stable
under Kazhdan--Lusztig tensoring with $\KL_\ka.$ 
\eth

We will need the following lemma.

\ble{ol} Every $\gt$-modules of central charge $\ka - h\spcheck$
such that
\begin{equation}
\label{prop}
V(k) \in \H, \; \mbox{respectively} \; \O'
\end{equation}
belongs to $\wh{\H}_{\ka, \ig},$ respectively $\wh{\O}'_{\ka, \ig}.$
\ele

\begin{proof} Fix a $\gt$-module with the property \eqref{prop}.
The action of $L_0$ on $V$ preserves $V(k).$ On each $V(k)/V(k-1),$ $L_0$ 
acts by $1/\ka \Om$ where $\Om$ is the Casimir operator of $\g.$
This implies that $L_0$ acts locally semisimply on $V.$ 
Denote the generalized eigenvalues of $L_0$ on $V(1)$
by $\xi_1,$ $\ldots,$ $\xi_n.$ Since $\g t$ maps $V^{\xi}$ to
$V^{\xi-1}$ and $V(k) \backslash V(1)$ isomorphically to $V(k-1)$
we have
\[
V^\xi = 0 \; \mbox{unless} \; \xi - \xi_i \in \Zset_{>0} \; 
{\mbox{for some}} \; i=1, \ldots,n.
\]
Then for any $\xi \in \Cset,$ $V^\xi \sub V(k)$ for some sufficiently 
large $k$ and therefore $V^\xi \in \H$ or $V^\xi \O'$ for the two 
cases in \eqref{prop}.
\end{proof}

\noindent
{\em{Proof of \thref{Fin}.}} Fix two $\gt$-modules $U \in \KL_\ka$ and 
$V \in \wh{\H}_\ka.$ Due to \leref{tensorle}
\[
T'(U,V)(k) \in \H, \; \forall k \in \Zset_{>0}.
\] 
\leref{ol} implies
\[
T'(U, V) \; \mbox{and} \; U \dot{\otimes} V \in \wh{\H}_{\ka, \ig}.
\]
We need to show that $U \dot{\otimes} V$ has finite length.

The proof of \thref{JKL} can be easily extended to show that 
\eqref{JKLisom} holds in this situation,
where in the lhs of \eqref{JKLisom} $\ja$ denotes the affine 
Jacquet functor from $\wh{\H}_{\ka, \ig}$ to $\wh{\O}'_{\ka ,\ig }.$

Since $V \in \wh{\H}_\ka$ from \prref{Jaind} we obtain 
$\ja(V) \in \wh{\O}^{'N}_\ka$ for some sufficiently large integer
$N.$ Then \leref{tensorle} implies that
\[
T'(U, D(\ja(V)))(k) \in \O^{'N}, \forall k \in \Zset_{>0}.
\]
Applying the characterization of the categories 
$\wh{\O}^{'N}_\ka$ from \prref{StrOO}, we get that 
$T'(U, D(\ja(V)))$ belongs to $\wh{\O}^{'N}_\ka.$ Then 
\eqref{JKLisom} gives that $\ja(U \dot{\otimes} V)$ belongs to 
$\wh{\O}'_\ka,$ as well, 
and in particular $\ja(U \dot{\otimes} V)$ has finite length.
Since the affine Jacquet functor is faithful we finally obtain that
$U \dot{\otimes} V$ has finite length.
\qed  
\sectionnew{Appendix: The Kazhdan--Lusztig fusion tensor product}
In this appendix we review the definition of the Kazhdan--Lusztig fusion 
tensor product \cite{KL1}.
Consider the Riemann sphere $\cp$ with three fixed distinct 
points $p_i,$ $i=0,1,2$ on it. Choose local coordinates (charts) at
each of them, i.e. isomorphisms $\ga_i: \cp \ra \cp$ such that
$\ga_i(p_i)= 0$ where the second copy of $\cp$ is equipped with 
a fixed coordinate function $\e$ vanishing at $0.$

Set $R=\Cset[\cp \backslash \{p_0, p_1, p_2\}]$
and denote by $\Ga$ the central extension
of the Lie algebra $\g \otimes R$
by
\begin{equation}
\label{Gaext}
[f_1 x_1,  f_2 x_2] := f_1 f_2 [x_1, x_2] +
\Res_{p_0} (f_2 d f_1) (x_1, x_2) \c,
\end{equation}
for $f_i \in R$ and $x_i \in \g.$ Here $(.,.)$ is invariant
bilinear form on $\g$ from section 2.3. There is a canonical
homomorphism 
\begin{equation}
\label{Gahom}
\Ga \ra \widehat{\g \oplus \g}, \quad
x f \mapsto (x \ev (\ga_1^\st)^{-1}(f), x \ev (\ga_2^\st)^{-1}(f)),
\; \c \mapsto -\c
\end{equation}
where $\ev(.)$ denotes the power series expansion of a rational 
function on $\cp$ at 0 in terms of the coordinate function $\e.$

Define 
\[
G_N = \Span \{ (f_1 x_1) \ldots (f_N x_N) \mid f_i 
\;\mbox{vanish at} \; p_0, \; x_i \in \g \} \sub U(\Ga).
\]

Let $f_0$ be a rational function on $\cp$ (unique up to a multiplication
by a nonzero complex number)
having only one (simple) zero
at $p_0$ and only one (simple) pole at $p_1.$ 
For instance when $\ga_0(p_1)$ is finite 
$f_0(\e) = a \ga^\st_0(\e)/(\ga^\st_0(\e)-\ga_0(p_1)),$ $a \neq 0.$ Set
\begin{equation}
\label{XN}
X_N = \Span \{ (f_0 x_1) \ldots (f_0 x_N) \mid x_i \in \g \} \sub U(\Ga).
\end{equation}
Clearly $X_N \sub G_N.$

Let $\f$ be a subalgebra of $\g$ which is reductive in $\g.$
Denote by $\H_{\g, \f}$ the category of finite length 
$\g$-modules which are $\f$-locally finite, semisimple, and admissible.
The last condition on a $\g$-module $M$ means that 
every finite dimensional $\f$-module appears with finite
multiplicity in  
$M|_{\f}.$ By $V \mt V^d$ we will denote the duality functor in 
$\H_{\g, \f}$ given by $M^d = (M^\st)^{U(\f)-fin}.$ By $D(.)$
we will denote the duality functor in $\Aff(\H_{\g, \f})_\ka$
given by $D(V)= (V^d)^\sharp(\infty),$ see \cite[Section 4]{Y}
for details.

Fix two smooth $\gt$ modules $U$ and $V$ of central charge
$\ka - h\spcheck$ such that
\[
U \in \KL_\ka \; \mbox{and} \;
V \in \Aff(\H_{\g, \k})_\ka.
\]

Using the homomorphism \eqref{Gahom} $U \otimes_\Cset V$ becomes 
a $\Ga$ module of central charge $-\ka + h\spcheck.$
Note that the restricted dual $(U \otimes V)^d$ 
is naturally a $\Ga$ submodule
of $(U \otimes V),$ having central charge
$\ka - h\spcheck.$ (The restricted dual is taken with respect to
the embedding $\k \hra \g \hra \Ga$ using constant functions
on $\cp.)$
Following \cite{KL1} define the following $\Ga$ submodule
of $(U \otimes V)^d$ 
\begin{equation}
\label{T}
T'(U, V) := \cup_{N \geq 1} T(U, V)\{N\} 
\end{equation}
where
\[
T(U, V)\{N\} := \Ann_{G_N} (U \otimes V)^d.
\] 
Eq. \eqref{T} indeed defines a $\Ga$ submodule of 
$(U \otimes V)^d$ since for any
$y \in \Ga$ and any integer $N$ there exists an integer $i$ such that
$G_{N+i} y \in \Ga G_N.$ In other words 
$T'(U,V)\{N\}$ is defined as 
\[
T'(U,V)\{N\} = \{ \eta \in (U \otimes V)^\st \mid \eta(G_N W) = 0, \; 
\dim U(\f) \eta < \infty \}.
\]

The spaces $T'(U, V)\{N\}$ inherit a natural $\g$-action,
using the embedding of $\g$ in $\Ga$ by constant functions.
Kostant's theorem \cite{Kos} that $\H_{\g, \k}$ is stable under tensoring 
with finite dimensional $\g$-modules and the following lemma of
Kazhdan and Lusztig imply that the spaces $T'(U, V)\{N\},$ 
equipped with this $\g$-action belong to $\H_{\g, \k}.$ 

\ble{tensorle}\cite[Proposition 7.4]{KL1} Assume that $V_i$ are two
strictly smooth $\gt$-modules of central charge $\ka-h\spcheck,$
generated by $V_i(N_i),$ respectively. Then
\[
U \otimes V = \sum_{k=0}^{N-1} X_k 
(U(N_1) \otimes V(N_2))+ 
G_N (U \otimes V)
\]
for all $N \in \Zset_{>0}.$
\ele

There is a canonical action of $\gt$ (``the copy attached to $p_0$'') 
on the space $T'(U,V),$ defined as follows. 
Let $\eta \in T(U,V)\{N\}.$ 
Fix $\om \in \Cset[\e, \e^{-1}],$ $x \in \g$ 
and choose $f \in R$ such that $f-\ga_0^\st(\om)$ has a zero of order at 
least $N$ at $p_0.$ Then
\begin{equation}
\label{tens}
(\om x) \eta := (f x) \eta
\end{equation}
correctly defines a structure of smooth $\gt$-module on $T'(U,V)$
of central charge $\ka - h\spcheck.$ Moreover 
\begin{equation}
\label{NN}
T'(U, V)(N) = T'(U, V)\{N\}.
\end{equation}

The Kazhdan--Lusztig fusion tensor product of the strictly 
smooth $\gt$-modules $U$ and $V$ is defined by
\begin{equation}
\label{KLproduct}
U \dot{\otimes} V = D T'(U,V) 
\end{equation}

In \cite{Y} we proved the following affine version of Kostant's
theorem that the categories $\H_{\g, \k}$ are stable under 
tensoring with finite dimensional $\g$-modules. 

\bth{Y} Assume that $\f$ is a subalgebra of $\g$ which is reductive in 
$\g$ and $\ka \notin \Rset_{\geq 0}.$ Then the categories 
$\Aff(\H_{\g, \f})_\ka$ are stable under Kazhdan--Lusztig 
tensoring with modules from the Kazhdan--Lusztig's category
$\KL_\ka.$ 
\eth

The hardest step in \thref{Y} is to show that $U \dot{\otimes} V$
is finitely generated (equivalently has finite length) for any 
$U \in \KL_\ka$ and $V \in \Aff(\H_{\g, \f})_\ka.$

The case $\f=\g$ recovers the finiteness result of Kazhdan and Lusztig
\cite{KL1} that $\KL_\ka$ is a monoidal category.

\end{document}